\documentclass{amsart}
\usepackage{amsmath,amscd,amssymb,amsthm,graphicx}

\def\p{p}
\def\en{e}
\def\Simp{\mathop{\rm Simp}\nolimits}

\def\ec{equivalence class}
\def\ep{\epsilon}
\def\even{{\mathop{\rm even}\nolimits}}
\def\odd{{\mathop{\rm odd}\nolimits}}

\def\gd{{\gamma, \delta}}
\def\gdp{{\gamma', \delta'}}

\def\lu{{\lambda, \mu}}
\def\lup{{\lambda', \mu'}}
\def\ul{{\mu, \lambda}}
\def\ulp{{\mu', \lambda'}}
\def\ll{{\lambda, \lambda}}

\def\EC{\mathop{\rm EC}\nolimits}

\def\CQ{\mathop{\rm CQ}\nolimits}
\def\pxi{\partial_\xi}
\def\Soo{{S^{1|1}}}
\def\Tan{\mathop{\rm Tan}\nolimits}
\def\VRoo{\Vec(\Roo)}

\def\cal{\mathcal}

\def\D{{\cal D}} 
 
\def\F{{\cal F}}

\def\I{{\cal I}}

\def\K{{\cal K}}

\def\frak{\mathfrak}

\def\dg{{\frak g}}

\def\dgo{{\frak o}}
\def\osp{\dgo\ds\dgp}
\def\dgp{{\frak p}}

\def\ds{{\frak s}}

\def\dU{{\frak U}}

\def\Bbb{\mathbb}

\def\bC{\Bbb C}

\def\bN{\Bbb N}

\def\bR{\Bbb R}

\def\bZ{\Bbb Z}

\def\h{\hat}

\def\lc{\lceil}
\def\lf{\lfloor}
\def\o{\overline}
\def\ot{\otimes}

\def\rc{\rceil}
\def\rf{\rfloor}
\def\t{\tilde}

\def\and{\mbox{\rm \ and\ }}

\def\Hom{\mathop{\rm Hom}\nolimits}

\def\lac{\mathop{\rm lac}\nolimits}

\def\mod{\mathop{\rm mod}\nolimits}

\def\oh{{\ts\frac{1}{2}}}

\def\px{\partial_x}

\def\Roo{{\bR^{1|1}}}

\def\Span{\mathop{\rm Span}\nolimits}

\def\SQ{\mathop{\rm SQ}\nolimits}

\def\Symb{\mathop{\rm Symb}\nolimits}

\def\ts{\textstyle}
\def\Vec{\mathop{\rm Vect}\nolimits}

\def\VR{\Vec(\Bbb R)}

\def\bgno{\bigbreak\noindent}

\def\cs{composition series}
\def\dog{differential operator}
\def\eg{{\em e.g.,\/}}

\def\ie{{\em i.e.,\/}}
\def\iff{if and only if}

\def\irr{irreducible}

\def\lsa{Lie superalgebra}

\def\meno{\medbreak\noindent}
\def\ncr{non-commutative residue}

\def\psdog{pseudo\dog}
\def\psidog{$\Psi$DO}
\def\r{representation}  

\def\sq{subquotient}
\def\tdm{tensor density module}

\def\th{\thinspace}

\def\uea{universal enveloping algebra}

\def\vf{vector field}

\newtheorem{lemma}{Lemma}[section]

\newtheorem{prop}[lemma]{Proposition}
\newtheorem{thm}[lemma]{Theorem}
\newtheorem{cor}[lemma]{Corollary}

\title[Supersymmetric subquotients]{Equivalence classes of subquotients of supersymmetric pseudodifferential operator modules}

\author{Charles H.\ Conley}
\address{Department of Mathematics \\University of North Texas \\Denton TX 76203, USA} 
\email{conley@unt.edu}

\thanks{The author was partially supported by Simons Foundation Collaboration Grant 207736.}

\begin{document}

\begin{abstract}
We study the \ec es of the non-resonant \sq s of spaces of \psdog s between \tdm s over the superline $\Roo$, as modules of the Lie superalgebra of contact \vf s.  There is a 2-parameter family of \sq s with any given Jordan-H\"older \cs.  We give a complete set of even equivalence invariants for \sq s of all lengths~$l$.  In the critical case $l = 6$, the even \ec es within each non-resonant 2-parameter family are specified by a pencil of conics.  In lengths $l \ge 7$ our invariants are not fully simplified: for $l = 7$ we expect that there are only finitely many equivalences other than conjugation, and for $l \ge 8$ we expect that conjugation is the only equivalence.  We prove this in lengths $l \ge 15$.  We also analyze certain lacunary \sq s.
\end{abstract}


\maketitle

\section{Introduction}  \label{Intro}

Given a manifold $M$, write $\Vec(M)$ for the Lie algebra of \vf s on $M$.  The natural $\Vec(M)$-module of functions on $M$ has a 1-parameter family of deformations, the {\em \tdm s\/} $\F_\lambda(M)$.  The space $\D_\lu(M)$ of \dog s from $\F_\lambda(M)$ to $\F_\mu(M)$ is a $\Vec(M)$-module under the adjoint action, and its {\em order filtration\/} $\D^k_\lu(M)$ is invariant.  Duval and Ovsienko initiated the investigation of the \ec es of the modules in this filtration in \cite{DO97}, treating the cases in which $\mu = \lambda$ and $k \le 2$.

The results of \cite{DO97} generalize in various directions.  In the multidimensional case, Lecomte, Mathonet, and Tousset extended them to all~$k$ \cite{LMT96}.  In the 1-dimensional case, Gargoubi considered the \ec es of $\D^k_\lu(\bR)$ for $\mu$ not necessarily equal to $\lambda$ \cite{Ga00}.  For $\mu$ equal to $\lambda$, Lecomte and Ovsienko treated \ec es of \sq\ modules $\D_\ll^k(M) / \D_\ll^{k-l}(M)$ of Jordan-H\"older length~$l$, and for $M = \bR$ they admitted \sq s $\Psi^k_\ll(\bR) / \Psi^{k-l}_\ll(\bR)$ of \psdog\ modules, $k \in \bC$ \cite{LO99}.  In general, the results depend only on the dimension of the manifold $M$.  The problem is not yet solved in complete generality even for $\bR$: the \ec es of the so-called {\em non-resonant\/} \sq s $\Psi_\lu^k(\bR) / \Psi_\lu^{k-l}(\bR)$ of lengths $l=6$ and~$7$ are not yet classified, and the resonant cases are unresolved for $l \ge 5$ \cite{CL13}.

In another direction, one may consider equivalence with respect to a subalgebra of $\Vec(M)$, for example, if $M$ has a contact structure, the Lie algebra $\K(M)$ of {\em contact \vf s.\/}  This problem was formulated for the superline $\Roo$ in \cite{GMO07}, and the \ec es of the spaces $\D^k_\lu(\Roo)$ regarded as modules of $\K(\Roo)$ were classified in \cite{BBB13}.

In the current article we study the \ec es of the \psdog\ \sq s $\Psi^{k,\p}_\lu(\Roo) / \Psi^{k-l/2,\, \p+l}_\lu(\Roo)$ under the action of $\K(\Roo)$, in the non-resonant cases.  Here the order $k$ is a complex parameter, $\p \in \bZ_2$ is a parity parameter, and $l$ is essentially the Jordan-H\"older length.  These \sq s exist because $\Psi^{k,\p}_\lu(\Roo)$ has a filtration ``twice as fine'' as the usual order filtration; see \cite{GMO07} and Lemma~\ref{D fine filt} below.  For example, the symbol modules $\Psi^{k,\p}_\lu(\Roo) / \Psi^{k-1, \p}_\lu(\Roo)$ are in general of length~2 as $\K(\Roo)$-modules.

The \cs\ of $\Psi^{k,\p}_\lu(\Roo) / \Psi^{k-l/2,\, \p+l}_\lu(\Roo)$ is determined by $\mu - \lambda - k$ and $l$, and its parity is determined by~$\p$.  Consequently there is a 2-parameter family of \sq s with any fixed \cs.  We find that in lengths $l \le 5$ almost all \sq s with a given \cs\ are equivalent, while in length~$6$ the even \ec es are determined by a single continuous invariant whose level curves form a pencil of conics.

In lengths $l \ge 7$ the even \ec es are determined by an essentially complete set of $2l - 11$ continuous invariants, each of whose level curves forms a pencil of conics.  Therefore it is natural to conjecture that in these lengths each \sq\ is evenly equivalent only to its conjugate, except possibly for a finite number of exceptional choices of length~$7$ composition series for which certain even \ec es are comprised of a finite number of conjugate pairs.  The higher the length, the more invariants one has, and so the easier this conjecture is to check.  A preliminary analysis using a software package shows that it holds in lengths $l \ge 15$, and we expect that it is not difficult for $l \ge 8$.  We would be interested to see complete results in length~$7$, where the necessary elimination theory is similar to that required in the unresolved length~$6$ and~$7$ cases over $\VR$ \cite{CL13}.

We work in the polynomial category: it is a general phenomenon that results in this area are the same in the smooth and polynomial categories.  Our strategy was outlined in Section~7.4 of \cite{Co09a} and Section~5 of \cite{Co09b} (but as we will see, the paragraph in \cite{Co09b} on the length~7 case is incorrect: there are three invariants rather than two).  Our central tool is Theorem~6.5 of \cite{Co09a}, which gives the action of $\K(\Roo)$ on the non-resonant \sq s of $\Psi^{k,\p}_\lu(\Roo)$ in terms of the {\em projective quantization.\/}

The organization of this article is roughly parallel to that of \cite{CL13}.  In Section~\ref{Dfns} we define non-resonant \sq s of \psdog\ modules, and in Section~\ref{ED} we recall conjugation, the duality resulting from the Berezin integral, and the super analog of the de Rham differential.  In Section~\ref{Lle6} we state our main result, the description of the \ec es of \sq s of lengths $l \le 6$, and in Section~\ref{Lge7} we give preliminary results on \sq s of length $l \ge 7$.  In Section~\ref{Lacunarity} we discuss the \sq\ question for {\em lacunary\/} \psdog\ modules.  Section~\ref{Proofs} contains the proofs, and in Section~\ref{Resonance} we conclude with some remarks on resonant \sq s.

\section{Definitions}  \label{Dfns}

The results in this section are given at greater length in \cite{GMO07} and \cite{Co09a}.  We denote the non-negative integers by $\bN$ and the positive integers by $\bZ^+$, and we use the Pochhammer symbol $(x)_r$ for the falling factorial:
\begin{equation*}
   (x)_r := x(x-1)\cdots (x-r+1), \qquad (x)_0 := 1.
\end{equation*}
For $x$ real, we write $\lf x \rf$ for the greatest integer less than or equal to~$x$, $\lc x \rc$ for the least integer greater than or equal to~$x$, and $\{x\}$ for the fractional part $x - \lf x \rf$.

Given a superspace $V$, we write $V_\even$ and $V_\odd$ for its even and odd parts, and $\ep_V: V \to V$ for the parity endomorphism $v \mapsto (-1)^{|v|} v$.  We write the parity-reversing functor $\Pi$ as a superscript:
\begin{equation*}
   V^\Pi := V, \qquad V^\Pi_\even := V_\odd, \qquad 
   V^\Pi_\odd := V_\even, \qquad \ep_{V^\Pi} = -\ep_V.
\end{equation*}
For $\p$ in $\bZ$ or $\bZ_2$, $V^{\p \Pi}$ denotes $V$ if $\p$ is even and $V^\Pi$ if $\p$ is odd.

Suppose that $V$ and $V'$ are superspaces carrying \r s $\pi$ and $\pi'$ of a \lsa\ $\dg$.  One has the usual notions of $\dg$-intertwining maps and $\dg$-equivalences from $V$ to $V'$, and they may be even, odd, or of mixed parity.  When the algebra $\dg$ is fixed by the context, we sometimes write $V \cong V'$ to indicate equivalent \r s.

We shall refer to the class of all \r s equivalent to $(\pi, V)$ by an even equivalence as its {\em even \ec.\/}  The even \ec\ of a \r\ is, of course, a subclass of its full \ec.  The space $V^\Pi$ carries the \r\ $\pi^\Pi$ defined by $\pi^\Pi(X) := \pi(X)$, and $\ep_V: V \to V^\Pi$ is an odd equivalence from $\pi$ to $\pi^\Pi$ (but the identity map is not).  Therefore odd equivalences from $V$ to $V'$ may be regarded as even equivalences from $V^\Pi$ to $V'$.

\subsection{Contact vector fields}

Henceforth we work exclusively over the superline, so we will drop the argument $\Roo$ and write simply $\K$, $\F_\lambda$, $\D_\lu$, $\Psi^{k,\p}_\lu$, and so on.  Let $\bC[x, \xi]$ be the polynomials on $\Roo$, with even coordinate $x$ and nilpotent odd coordinate $\xi$ (so $\xi^2 = 0$).  The \lsa\ $\VRoo$ of polynomial \vf s on $\Roo$ is $\Span_{\bC[x, \xi]} \{\px, \pxi\}$.

The standard {\em contact 1-form\/} on $\Roo$ is
\begin{equation*}
   \alpha := dx + \xi d\xi.
\end{equation*}
The kernel of this form is a completely non-integrable distribution $\Tan$ on $\Roo$, and the contact subalgebra $\K$ of $\VRoo$ is defined to be the stabilizer of $\Tan$.  To be explicit, set
\begin{equation*}
   D:= \pxi + \xi \px, \qquad \o D := \pxi - \xi \px,
\end{equation*}
and let $X: \bC[x, \xi] \to \VRoo$ be the even injection
\begin{equation*}
   X\bigl(F(x, \xi)\bigr) := F \px + \oh D(F) \o D.
\end{equation*}
We will often write $X_F$ for $X(F)$.  One obtains
\begin{equation*}
   \K = X \bigl(\bC[x, \xi] \bigr), \qquad \Tan = \bC[x, \xi] \o D, \qquad
   \VRoo = \K \oplus \Tan,
\end{equation*}
the direct sum being $\K$-invariant (but not $\bC[x, \xi]$-invariant).

\subsection{Tensor densities}

The {\em \tdm s\/} $\{ \F_\lambda: \lambda \in \bC \}$ of $\K$ comprise a 1-parameter family of deformations of the natural module $\bC[x, \xi]$.  We write $\F_\lambda$ as
\begin{equation*}
   \F_\lambda := \alpha^\lambda\, \bC[x, \xi],
\end{equation*}
where $\alpha^\lambda$ is defined to be even so $\F_{\lambda, \even} = \alpha^\lambda \bC[x]$ and $\F_{\lambda, \odd} = \alpha^\lambda \xi \bC[x]$.  The Lie action $L_\lambda$ of $\K$ on $\F_\lambda$ is
$$ L_\lambda (X_F) (\alpha^\lambda G) :=
   \alpha^\lambda \bigl( X_F(G) + \lambda \px(F) G \bigr). $$
It is well-known that $\F_\lambda$ is \irr\ under $\K$ unless $\lambda=0$; $\F_0$ contains the trivial submodule $\bC$, and $\F_0 / \bC$ is evenly equivalent to $\F_{1/2}^\Pi$.  Moreover, no two of the $\F_\lambda$ are equivalent.

\subsection{Differential operators}

Let $\D_\lu$ be the space $\alpha^{\mu-\lambda} \bC[x, \xi, \px, \pxi]$ of \dog s from $\F_\lambda$ to $\F_\mu$.  Given any element $\tau$ of $\bC[x, \xi, \px, \pxi]$, the action of $\alpha^{\mu-\lambda} \tau$ on $\F_\lambda$ is defined by $(\alpha^{\mu-\lambda} \tau) (\alpha^\lambda G) := \alpha^\mu \tau(G)$.  It is conventional to write $\delta$ for the difference $\mu - \lambda$:
\begin{equation*}
   \delta(\lu) := \mu - \lambda.
\end{equation*}

Using the fact that $\o D^2 = -\px$, one finds that $\{ \alpha^\delta \o D^j: j \in \bN \}$ is a basis of $\D_\lu$ over $\bC[x,\xi]$: given any element $T$ of $\D_\lu$, there are unique elements $T_j$ of $\bC[x, \xi]$, all but finitely many equal to zero, such that $T = \alpha^\delta\, \sum_{j \in \bN} T_j \o D^j$.

The Lie action $L_\lu$ of $\K$ on $\D_\lu$ is
\begin{equation} \label{Llu}
   L_\lu (X_F) (T) :=
   L_\mu (X_F) \circ T - T \circ L_\lambda (X_F).
\end{equation}
The formula for this action in terms of the coefficients $T_j$ is complicated and will not be needed in this work.

For $k \in \oh\bN$, define
\begin{equation*}
   \D_\lu^k := \Span_{\bC[x, \xi]} \bigl\{ \alpha^\delta \o D^j: 0 \le j \le 2k \bigr\}.
\end{equation*}
The following lemma is proven in \cite{GMO07}.  Part~(i) follows from the facts that $\K_\odd$ generates $\K$, $X_{\xi f(x)} = \oh f D$, and $[D, \o D] = 0$.  Part~(ii) follows from $\o D^2 = -\px$, and Part~(iii) is a computation using~(\ref{Llu}).

\begin{lemma} \label{D fine filt}
\begin{enumerate}
\item[(i)]
$\D^k_\lu$ is $\K$-invariant for all $k \in \oh\bN$.
\item[(ii)]
For $k \in \bN$, $\D_\lu^k$ is the space of \dog s of order~$\le k$.
\item[(iii)]
There is an even $\K$-equivalence from $\D^k_\lu / \D^{k - 1/2}_\lu$ to $\F_{\delta - k}^{2k\Pi}$, defined by
\begin{equation*}
   \alpha^\delta\, \sum_{j=0}^{2k} T_j \o D^j \mapsto \alpha^{\delta - k} T_{2k}.
\end{equation*}
\end{enumerate}
\end{lemma}

\subsection{Pseudodifferential operators}

For $z \in \bC$ and $\p \in \bZ_2$, define formal symbols $\o D^z_\p$ of parity $\p$:
\begin{equation*}
   \o D^z_0 := e^{i \pi z/2}\, \px^{z/2}, \qquad
   \o D^z_1 := e^{i \pi (z-1)/2}\, \px^{(z-1)/2} \o D.
\end{equation*}
Note that $\o D^{z'}_{\p'} \circ \o D^z_\p = \o D^{z'+z}_{\p'+\p}$, and for $j \in \bN$, $\o D^j_{j \mod 2}$ is simply $\o D^j$.

\meno {\bf Definition.}
For $k \in \bC$ and $\p \in \bZ_2$, the space $\Psi^{k,\p}_\lu$ of {\em \psdog s\/} (\psidog s) of order~$\le (k,\p)$ from $\F_\lambda$ to $\F_\mu$ consists of formal series:
\begin{equation*}
   \Psi^{k,\p}_\lu := \Bigl\{ \alpha^\delta \sum_{j \in \bN}
   T_{2k-j} \o D^{2k-j}_{ \p + (j \mod 2)}:\, T_{2k-j} \in \bC[x, \xi] \Bigr\}.
\end{equation*}
Observe that $\Psi^{k,\p}_\lu \subset \Psi^{k + 1/2,\, \p+1}_\lu \subset \Psi^{k+1,\, \p}_\lu \subset \cdots$.  We name the nested union: 
\begin{equation*}
   \Psi^{\bN+k,\, \p}_\lu := \bigcup_{i=0}^\infty \Psi^{k+i,\, \p}_\lu.
\end{equation*}

\medbreak
Although \psidog s are not actually operators on \tdm s, composition of \dog s extends to a composition map from $\Psi^{k', \p'}_{\mu, \nu} \ot \Psi^{k, \p}_\lu$ to $\Psi^{k+k',\, \p+\p'}_{\lambda, \nu}$ via $[\o D, \px^z] = 0$ and the {\em generalized Leibniz rule\/}
\begin{equation*}
   \px^z \circ F(x, \xi) := \sum_{j \in \bN} 
   {\ts {z \choose j}} \px^j(F) \circ \px^{z-j}.
\end{equation*}
Thus~(\ref{Llu}) extends to define an action $L_\lu$ of $\K$ on $\Psi^{\bN+k,\, \p}_\lu$.  The proof of Lemma~\ref{D fine filt} also proves the following lemma, which is Lemma~6.2 in \cite{Co09a}.

\begin{lemma} \label{Psi fine filt}
\begin{enumerate}
\item[(i)]
$\Psi^{k,\p}_\lu$ is $\K$-invariant for all $k \in \bC$ and $\p \in \bZ_2$.
\item[(ii)]
There is an even $\K$-equivalence from $\Psi^{k,\p}_\lu / \Psi^{k - 1/2,\, \p+1}_\lu$ to $\F_{\delta - k}^{\p\Pi}$, defined by
\begin{equation*}
   \alpha^\delta\, \sum_{j \in \bN} T_{2k-j} \o D^{2k-j}_{\p+j} \mapsto \alpha^{\delta - k} T_{2k}.
\end{equation*}
\end{enumerate}
\end{lemma}

\subsection{Subquotients} \label{Sbqnts}

For $\lambda$, $\mu$, and $k$ in $\bC$, $\p$ in $\bZ_2$, and $l$ in $\bZ^+$, we define the $\K$-module $\SQ^{k,\p,l}_\lu$ by
\begin{equation*}
   \SQ^{k,\p,l}_\lu := \Psi^{k,\p}_\lu\, \big/\, \Psi^{k-l/2,\, \p+l}_\lu.
\end{equation*}
The topic of this paper is the description of the equivalence classes of the collection of all of these \sq\ modules.  By Lemma~\ref{Psi fine filt}, $\SQ^{k,\p,l}_\lu$ has a filtration with successive \sq s
\begin{equation*}
   \bigl\{ \F_{\delta-k}^{\p\Pi},\, \F_{\delta-k+1/2}^{(\p+1)\Pi},\, \ldots,\,
              \F_{\delta-k+(l-1)/2}^{(\p+l-1)\Pi} \bigr\}.
\end{equation*}
We shall refer to $\SQ^{k,\p,l}_\lu$ as having Jordan-H\"older length~$l$, even though this is slightly inaccurate if $0 \in \{\delta-k, \ldots, \delta-k+l-1\}$ because $\F_0$ is itself of length~2 rather than \irr.  Nevertheless, the following lemma is immediate.

\begin{lemma} \label{comp series}
If\/ $\SQ^{k,\p,l}_\lu \cong \SQ^{k',\p',l'}_\lup$, then $l = l'$, $\delta - k = \delta' - k'$, and the parity of the equivalence is $\p - \p'$.  There are no equivalences of mixed parity.
\end{lemma}

In light of this necessary condition for equivalence, our topic reduces to the following ``equivalence question''.  Define
\begin{equation*}
   n(k, \delta) := \delta - k.
\end{equation*}

\meno {\bf Question.}
For fixed $n \in \bC$ and $l \in \bZ^+$, what are the $\K$-\ec es of the set
\begin{equation*}
   \bigl\{ \SQ^{\delta-n,\, \p,\, l}_\lu:\, \lambda,\, \mu \in \bC,\, \p \in \bZ_2 \bigr\}
\end{equation*}
of length~$l$ \sq s with \cs\ $\{ \F_n^{\p\Pi},\, \F_{n + 1/2}^{(\p+1)\Pi},\, \ldots,\, \F_{n + (l-1)/2}^{(\p+l-1)\Pi}\}$?

\meno{\bf Remark.}
The equivalence question includes the question of the \ec es of the \dog\ modules addressed in \cite{BBB13}, because for $k \in \oh\bN$, $\SQ^{k,\, 2k,\, 2k+1}_\lu$ is simply $\D^k_\lu$.  We shall recover the results of \cite{BBB13} in the non-resonant cases (see below) in Section~\ref{BBB Results}.

\subsection{Resonance}

In this article we consider only {\em non-resonant\/} equivalence classes.  {\em Resonance\/} is the failure of complete reducibility under the action of the {\em projective subalgebra\/} of $\K$.  This subalgebra is isomorphic to $\osp_{1|2}$ and is
\begin{equation*}
   \ds := \Span_\bC \bigl\{ X_1,\, X_\xi,\, X_x,\, X_{\xi x},\, X_{x^2} \bigr\}.
\end{equation*}

The Casimir operator $Q_\ds$ of $\ds$ spans the space of quadratic central elements of the \uea\ $\dU(\ds)$ with no constant term.  It is defined by
\begin{equation*}
   Q_\ds := T_\ds^2 - {\ts\frac{1}{16}}, \qquad
   T_\ds := X_x - 4 X_{\xi x} X_\xi - {\ts\frac{1}{4}}.
\end{equation*}
$T_\ds$ acts on $\F_\nu$ by $(\nu - \frac{1}{4}) \ep_{\F_\nu}$, and so $Q_\ds$ acts on both $\F_\nu$ and $\F_\nu^\Pi$ by
\begin{equation*}
   L_\nu(Q_\ds) = L_\nu^\Pi(Q_\ds) = \nu^2 - \oh \nu.
\end{equation*}

The module $\SQ^{\delta-n,\, \p,\, l}_\lu$ cannot be resonant unless its composition series has repeated Casimir eigenvalues.  Since the values of $\nu^2 - \oh \nu$ are symmetric around $\nu = \frac{1}{4}$, this occurs \iff\ $(n+i/2) + (n+j/2)$ is $\oh$ for some $0 \le i < j \le l-1$.  In fact a further condition is necessary for resonance: $2n$ must be integral.  If $n$ meets both of these conditions, then normally $\SQ^{\delta-n,\, \p,\, l}_\lu$ is not completely reducible under $\ds$.  More precisely, adapting the proof of Lemma~6.6 of \cite{GMO07}, we see that it is completely $\ds$-reducible \iff\ either $p + 2n$ is even and $\mu = \lambda$, or $p + 2n$ is odd and $\mu + \lambda = \oh$ (the self-adjoint case).  Therefore we make the following definition.

\meno{\bf Definition.}
The $l-1$ $n$-values $0, -\oh, -1, -\frac{3}{2}, \ldots, 1 - \oh l$ are {\em resonant with respect to~$l$.\/}  For these values of $n$, the subquotient $\SQ^{\delta-n,\, \p,\, l}_\lu$ is called {\em resonant.\/}

\section{Equivalences and dualities} \label{ED}

In this section we describe certain properties of equivalence classes that hold in all lengths.  These properties will permit us to answer the equivalence question more transparently in terms of $l$, $\p$, $\delta(\lu)$, and the following new parameters:
\begin{equation*}
   \gamma(\lu) := 3(\lambda + \mu -\oh)^2, \qquad
   N_l(n) := n + {\ts\frac{1}{4}} l - \oh.
\end{equation*}

\subsection{Conjugation of \psdog s}

Conjugation is the map $T \mapsto T^*$ from $\Psi^{k,\p}_\lu$ to $\Psi^{k,\p}_{-\mu + 1/2,\, -\lambda + 1/2}$ defined by
\begin{equation*}
   (\alpha^\delta G \o D_\p^z)^* :=
   e^{i\pi (z + \p)/2} (-1)^{\p |G|} \alpha^\delta \o D_\p^z G.
\end{equation*}
Note that conjugating twice acts on $\Psi^{k,\p}_\lu$ as the scalar map $e^{i\pi (2k+\p)}$.  The following lemma is well-known: see, \eg\ Section~4.2 of \cite{GMO07} or Proposition~6.7 of \cite{Co09a}.

\begin{lemma} \label{conj}
Conjugation is an even $\K$-equivalence from $\Psi^{k,\p}_\lu$ to $\Psi^{k,\p}_{-\mu + 1/2,\, -\lambda + 1/2}$.  In particular, there is an even equivalence
\begin{equation*}
   \SQ^{\delta-n,\, \p,\, l}_\lu \cong \SQ^{\delta-n,\, \p,\, l}_{ -\mu + 1/2,\, -\lambda + 1/2}.
\end{equation*}
\end{lemma}

\subsection{Dualities}

Let us temporarily pass to the category of $\K(\Soo)$-modules.  Algebraically, this is accomplished by simply adjoining $x^{-1}$ to all the objects under consideration.  We will denote the resulting extensions with a hat: $\h\K$, $\h\F_\lambda$, $\h\Psi^{k,\p}_\lu$, $\h\SQ^{k,\p,l}_\lu$, and so on.  The following results are collected in Section~6.2 of \cite{Co09a}.

There is, up to a scalar, a unique non-zero $\h\K$-map from $\h\F_{1/2}$ to the trivial even module $\bC^{1|0}$, the {\em Berezin integral.\/}  It is odd and can be used to define a non-degenerate symmetric odd $\h\K$-invariant pairing between $\h\F_\lambda$ and $\h\F_{-\lambda + 1/2}$.  The conjugation equivalence above is the adjoint map with respect to this pairing.

The Berezin integral also yields the {\em super \ncr,\/} a non-zero even $\h\K$-map from $\h\Psi^{\bN,\, 0}_{\lambda, 0}$ to $\bC^{1|0}$.  Combining the super \ncr\ with composition of \psdog s gives a non-degenerate supersymmetric even $\h\K$-invariant pairing between $\h\Psi^{\bN + k,\, \p}_\lu$ and $\h\Psi^{\bN - k,\, \p}_\ul$, the {\em super Adler trace.\/}

Using the fact that the dual of $\h\F_\lambda$ is $\h\F_{-\lambda + 1/2}^\Pi$, it can be shown that the super Adler trace drops to \sq s so as to give the following lemma.  The important point for us is the subsequent corollary, which applies to $\K$- rather than $\h\K$-\sq s.

\begin{lemma} \label{circle dual}
The $\h\K$-modules\/ $\h\SQ^{k,\p,l}_\lu$ and\/ $\h\SQ^{-k-1-l/2,\, \p+l,\, l}_\ul$ are dual.
\end{lemma}

\begin{cor} \label{line dual}
There is a $\K$-equivalence
\begin{equation*}
   \SQ^{\mu - \lambda - n,\, \p,\, l}_\lu\ \cong\
   \SQ^{\mu' - \lambda' - n,\, \p',\, l}_\lup
\end{equation*}
\iff\ there is a $\K$-equivalence
\begin{equation*}
   \SQ^{\lambda - \mu + n - 1 + l/2,\, \p + l,\, l}_\ul\ \cong\
   \SQ^{\lambda' - \mu' + n - 1 + l/2,\, \p' + l,\, l}_\ulp.
\end{equation*}
\end{cor}

\meno {\em Proof.\/}
First use Lemma~\ref{circle dual} to prove the analogous statement over $\Soo$.  Then deduce that any $\h\K$-isomorphism from $\h\SQ^{\mu - \lambda - n,\, \p,\, l}_\lu$ to $\h\SQ^{\mu' - \lambda' - n,\, \p',\, l}_\lup$ must carry $\SQ^{\mu - \lambda - n,\, \p,\, l}_\lu$ to $\SQ^{\mu' - \lambda' - n,\, \p',\, l}_\lup$, because these are the subspaces on which $X_1 = \px$ acts locally nilpotently.  $\Box$

\medbreak
We remark that for non-resonant modules, Corollary~\ref{line dual} follows directly from Theorem~6.5 and Proposition~6.12 of \cite{Co09a}.

For fixed $l$, $n$, and $\p$, Lemma~\ref{conj} implies that the even \ec\ of $\SQ^{\delta-n,\, \p,\, l}_\lu$ depends only on $(\gd)$, rather than on $(\lu)$.  By Corollary~\ref{line dual}, the equations defining this even \ec\ within the set of all \sq s are symmetric under $(N_l,\, \p,\, \gamma,\, \delta) \mapsto (-N_l,\, \p+l,\, \gamma,\, -\delta)$.  In particular, in even lengths they are symmetric under $(N_l,\, \delta) \mapsto (-N_l, -\delta)$, and in odd lengths this transformation exchanges the equivalence conditions in the two parities.  Therefore, as stated above, we will give the equations in terms of $N_l$, $\gamma$, and $\delta$.

Note that the set of resonant values of $N_l$ is $\{-\frac{1}{4} l + \oh,\, -\frac{1}{4} l + 1,\, \ldots,\, \frac{1}{4} l - \oh \}$.  Its symmetry around zero is a consequence of Lemma~\ref{circle dual}.  

Keep in mind that $(\gd)$ specifies a conjugate pair of values of $(\lu)$ rather than a single value, except in the self-adjoint cases where $\gamma = 0$.  In fact, some of our formulas involve $\gamma^{1/2}$.  Although the statements of our main results are independent of the choice of the sign of the square root, for concreteness we specify
\begin{equation*}
   \gamma^{1/2}(\lu) := \sqrt{3}\, (\lambda + \mu - \oh).
\end{equation*}

Henceforth we will always use the notation
\begin{equation*}
   (\gdp) := \bigl( \gamma(\lup), \delta(\lup) \bigr).
\end{equation*}
We make the following definition in order to be able to regard the even \ec\ of $\SQ^{\delta-n,\, \p,\, l}_\lu$ as a subset of the $(\gd)$-plane.

\meno {\bf Definition.}
$\EC^{\p,\, l}_n (\gd) := \bigl\{ (\gamma',\delta') \in \bC^2:\ 
\SQ^{\delta'-n,\, \p,\, l}_\lup \cong \SQ^{\delta-n,\, \p,\, l}_\lu \bigr\}$.

\subsection{The first super Bol operator}

Non-scalar maps from $\F_\lambda$ to $\F_\mu$ intertwining the actions of the projective subalgebra $\ds$ exist exclusively in the half-integral self-adjoint cases, where $\gamma = 0$ and $\delta \in \oh + \bN$.  For each such $\delta$, up to a scalar the unique such map is the {\em super Bol operator,\/} $\alpha^\delta \o D^{2\delta}$.

The only $\K$-invariant super Bol operator is the first one, at $\delta = \oh$.  It is the unique non-scalar $\K$-intertwining map between any two \tdm s and may be thought of as an analog of the de Rham differential, so we denote it by $\o d$:
\begin{equation*}
   \o d := \alpha^{1/2} \o D: \F_0 \to \F_{1/2}.
\end{equation*}
It gives rise to odd equivalences between certain \sq s of arbitrary length, which we now describe.

Write $L_{\o d}$ and $R_{\o d}$ for left and super-right composition with $\o d$, respectively:
\begin{equation*} \begin{array}{rl}
   L_{\o d}: \Psi^{k,\, \p}_{\lambda,\, 0} \to \Psi^{k+1,\, \p+1}_{\lambda,\, 1/2}, &
   \quad T \mapsto \o d \circ T, \\[6pt]
   R_{\o d}: \Psi^{k,\, \p}_{1/2,\, \mu} \to \Psi^{k+1,\, \p+1}_{0,\, \mu}, &
   \quad T \mapsto (-1)^{|T|} T \circ \o d.
\end{array} \end{equation*}
These maps are both odd $\K$-equivalences, which induce odd $\K$-equivalences
\begin{equation*}
   L_{\o d}: \SQ^{-\lambda - n,\, \p,\, l}_{\lambda,\, 0}
   \to \SQ^{-\lambda - n + 1/2,\, \p+1,\, l}_{\lambda,\, 1/2}, \quad
   R_{\o d}: \SQ^{\mu - n - 1/2,\, \p,\, l}_{1/2,\, \mu}
   \to \SQ^{\mu - n,\, \p+1,\, l}_{0,\, \mu}.
\end{equation*}
Observe that the two cases form a conjugate pair, so in $(\gd)$-coordinates they appear as a single case.  Thus we have:

\begin{lemma} \label{Bol}
For all $l$, $\p$, $n$, and $\nu$, there are odd equivalences from the elements of\/ $\EC^{\p,\, l}_n \bigl( 3(\nu + \oh)^2, \nu \bigr)$ to the elements of\/ $\EC^{\p+1,\, l}_n \bigl( 3\nu^2, \nu + \oh \bigr)$.
\end{lemma}

Now observe that for $k \in \oh \bN$ and $l \ge 2k+2$ we have the canonical $\K$-splittings
\begin{equation} \label{Dsplitting}
   \Psi^{\bN,\, 0}_\lu = \D_\lu \oplus \Psi^{-1/2,\, 1}_\lu, \qquad
   \SQ^{k,\, 2k,\, l}_\lu = \D^k_\lu \oplus \SQ^{-1/2,\, 1,\, l - 2k-1}_\lu.
\end{equation}
The interplay between this equation and the maps $L_d$ and $R_d$ gives the following lemma, which explains the exceptional splittings of certain \sq\ modules.

\begin{lemma} \label{DBol}
For all $\lambda$ and $\mu$, we have the\/ $\K$-splittings
\begin{equation*} \begin{array}{rclrcl}
   \D_{\lambda,\, 1/2} &=& \D^0_{\lambda,\, 1/2}
   \oplus L_{\o d} (\D_{\lambda,\, 0}), &
   \Psi^{-1/2,\, 1}_{\lambda,\, 0} &=& \Psi^{-1,\, 0}_{\lambda,\, 0}
   \oplus L_{\o d}^{-1} (\D^0_{\lambda,\, 1/2}), \\[6pt]
   \D_{0,\, \mu} &=& \D^0_{0,\, \mu}
   \oplus R_{\o d} (\D_{1/2,\, \mu}), &
   \Psi^{-1/2,\, 1}_{1/2,\, \mu} &=& \Psi^{-1,\, 0}_{1/2,\, \mu}
   \oplus R_{\o d}^{-1} (\D^0_{0,\, \mu}), \\[6pt]
   \D_{0,\, 1/2} &=& \D^{1/2}_{0,\, 1/2}
   \oplus L_{\o d} R_{\o d} (\D_{1/2,\, 0}), &
   \Psi^{-1/2,\, 1}_{1/2,\, 0} &=& \Psi^{-3/2,\, 1}_{1/2,\, 0}
   \oplus L_{\o d}^{-1} R_{\o d}^{-1} (\D^{1/2}_{0,\, 1/2}).
\end{array} \end{equation*}
\end{lemma}

\section{Lengths $l \le 6$}  \label{Lle6}

In this section we state our main result, the answer to the equivalence question in the non-resonant cases of length $l \le 6$.  The proof will be given in Section~\ref{Proofs}.

Recall that $n$, $\p$, and $l$ are invariants of the even \ec\ of $\SQ^{\delta - n,\, \p,\, l}_\lu$, which is represented by the subset $\EC^{\p,\, l}_n(\gd)$ of $\bC^2$.  By Lemma~\ref{comp series} there are two possibilities for the full \ec\ of the \sq: if there is an odd equivalence from $\SQ^{\delta - n,\, \p,\, l}_\lu$ to some \sq\ $\SQ^{\delta' - n,\, \p+1,\, l}_\lup$, then the full \ec\ is the union of the sets of \sq s represented by $\EC^{\p,\, l}_n(\gd)$ and $\EC^{\p+1,\, l}_n(\gdp)$, while if there is no such odd equivalence, then $\EC^{\p,\, l}_n(\gd)$ represents the full \ec.

\subsection{Lengths $l \le 5$} \label{Lle5}

In lengths $l \le 3$ there is nothing to prove: it follows from Theorem~6.5 and Corollary~6.6 of \cite{Co09a} that the non-resonant \sq s all split fully.  In length~$2$ the resonant value is $n = N_2 = 0$; for all other $n$ we have
\begin{equation*}
   \SQ^{\delta - n,\, \p,\, 2}_\lu \cong
   \F_n^{\p\Pi} \oplus \F_{n + 1/2}^{(\p+1)\Pi}.
\end{equation*}
In length~$3$ the resonant values are $n = 0$ and $-\oh$, \ie\ $N_3 = \pm \frac{1}{4}$.  For all other~$n$,
\begin{equation*}
   \SQ^{\delta - n,\, \p,\, 3}_\lu \cong
   \F_n^{\p\Pi} \oplus \F_{n + 1/2}^{(\p+1)\Pi} \oplus \F_{n+1}^{\p\Pi}.
\end{equation*}
Therefore in these lengths all non-resonant \sq s with a given \cs\ are equivalent: $\EC^{\p,\, l}_n(\gd) = \bC^2$ for $l \le 3$, and there is an odd equivalence between the two even \ec es of opposite parities. 

In order to state the results in higher lengths we make the following definitions:
\begin{equation} \label{Bij} \begin{array}{rcl}
   B^0_{m + 3/2,\, m}(\gd) &:=&
   \gamma^{1/2}, \\[6pt]
   B^1_{m + 3/2,\, m}(\gd) &:=&
   \gamma - \bigl[ 3(m + \oh) \delta + {\ts\frac{3}{4}} \bigr], \\[6pt]
   B^0_{m + 2,\, m}(\gd) &:=&
   \gamma - \bigl[ (m + {\ts\frac{3}{2}}) (2\delta + m + \oh) \bigr], \\[6pt]
   B^1_{m + 2,\, m}(\gd) &:=&
   \gamma - \bigl[ m(2\delta + m + 1) \bigr], \\[6pt]
   B^0_{m + 5/2,\, m}(\gd) &:=&
   \gamma - \bigl[ (m + 1) \delta + {\ts\frac{3}{4}} \bigr], \\[6pt]
   B^1_{m + 5/2,\, m}(\gd) &:=&
   \gamma^{3/2} - \gamma^{1/2}
   \bigl[ 4(m + 1) \delta - (m + 1)^2 + 3 \bigr].
\end{array} \end{equation}

\bgno {\bf Definition.}
Two \sq s $\SQ^{\delta-n,\, \p,\, l}_\lu$ and $\SQ^{\delta'-n,\, \p,\, l}_\lup$ are said to {\em induce simultaneous vanishing of the functions\/} $f_1(\gd), \ldots, f_r(\gd)$ if for all~$s$, $f_s(\gd)$ and $f_s(\gdp)$ are either both zero or both non-zero.

\medbreak
In length~$4$ the set of resonant values of $n$ is $\{-1, -\oh, 0 \}$, so that of $N_4$ is $\{0, \pm\oh \}$.  The description of $\EC^{\p, 4}_n(\gd)$ involves $B^\p_{n+3/2,\, n}$, which we rewrite in terms of $N_4$:
\begin{equation*}
   B^0_{n+3/2,\, n} = \gamma^{1/2}, \qquad
   B^1_{n+3/2,\, n} = \gamma - \bigl[ 3N_4 \delta + {\ts\frac{3}{4}} \bigr].
\end{equation*}

\begin{prop} \label{L4 prop}
Fix $n$ non-resonant with respect to $l=4$.
\begin{enumerate}
\item[(i)]
$\SQ^{\delta-n,\, 0,\, 4}_\lu \cong \SQ^{\delta'-n,\, 0,\, 4}_\lup$ \iff\ they induce simultaneous vanishing of 
\begin{equation} \label{L4p0}
   (\delta - N_4 + \oh)_2\, B^0_{n + 3/2,\, n}.
\end{equation}
The \ec\ where~(\ref{L4p0}) vanishes is $\F_n \oplus \F^\Pi_{n+1/2} \oplus \F_{n+1} \oplus \F^\Pi_{n+3/2}$.
\smallbreak\item[(ii)]
$\SQ^{\delta-n,\, 1,\, 4}_\lu \cong \SQ^{\delta'-n,\, 1,\, 4}_\lup$ \iff\ they induce simultaneous vanishing of 
\begin{equation} \label{L4p1}
   (\delta - N_4)\, B^1_{n + 3/2,\, n}.
\end{equation}
The \ec\ where~(\ref{L4p1}) vanishes is $\F^\Pi_n \oplus \F_{n+1/2} \oplus \F^\Pi_{n+1} \oplus \F_{n+3/2}$.
\smallbreak\item[(iii)]
$\SQ^{\delta-n,\, 0,\, 4}_\lu \cong \SQ^{\delta'-n,\, 1,\, 4}_\lup$ \iff~(\ref{L4p0}) at $(\gd)$ and~(\ref{L4p1}) at $(\gdp)$ are either both zero or both non-zero.
\end{enumerate}
\end{prop}

Note that Proposition~\ref{L4 prop} exhibits the $(N_4, \delta)$-symmetry discussed below Corollary~\ref{line dual}: $l$ is even, and~(\ref{L4p0}) and~(\ref{L4p1}) are even and odd in $(N_4, \delta)$, respectively.  Also, the coefficients $(\delta - N_4 + \oh)_2$ and $(\delta - N_4)$ of $B^\p_{n+3/2,\, n}$ in~(\ref{L4p0}) and~(\ref{L4p1}) are zero precisely when the \sq s are split by~(\ref{Dsplitting}).

In length~$5$ the set of resonant values of $n$ is $\{-\frac{3}{2}, -1, -\oh, 0 \}$, so that of $N_5$ is $\{\pm\frac{1}{4}, \pm\frac{3}{4} \}$.  The description of $\EC^{\p, 5}_n(\gd)$ involves $B^\p_{n+3/2,\, n}$, $B^{\p+1}_{n+2,\, n+1/2}$, and $B^\p_{n+2,\, n}$, which we now give in terms of $N_5$.  For $\p = 0$, $B^0_{n+3/2,\, n} = \gamma^{1/2}$,
\begin{equation*}
   B^1_{n+2,\, n+1/2} = \gamma - \bigl[ 3(N_5 + {\ts\frac{1}{4}}) \delta
   + {\ts\frac{3}{4}} \bigr], \quad
   B^0_{n+2,\, n} = \gamma - \bigl[ (N_5 + {\ts\frac{3}{4}})
   (2\delta + N_5 - {\ts\frac{1}{4}}) \bigr],
\end{equation*}
and for $\p = 1$, $B^0_{n+2,\, n+1/2} = \gamma^{1/2}$,
\begin{equation*}
   B^1_{n+3/2,\, n} = \gamma - \bigl[ 3(N_5 - {\ts\frac{1}{4}}) \delta
   + {\ts\frac{3}{4}} \bigr], \quad
   B^1_{n+2,\, n} = \gamma - \bigl[ (N_5 - {\ts\frac{3}{4}})
   (2\delta + N_5 + {\ts\frac{1}{4}}) \bigr].
\end{equation*}

\begin{prop} \label{L5 prop}
Fix $n$ non-resonant with respect to $l=5$.
\begin{enumerate}
\item[(i)]
$\SQ^{\delta-n,\, 0,\, 5}_\lu \cong \SQ^{\delta'-n,\, 0,\, 5}_\lup$ \iff\ they induce simultaneous vanishing of 
\begin{equation} \label{L5p0}
   (\delta - N_5 + {\ts\frac{3}{4}})_2\, B^0_{n + 3/2,\, n}, \quad
   (\delta - N_5 - {\ts\frac{1}{4}})\, B^1_{n + 2,\, n + 1/2}, \quad
   (\delta - N_5 + {\ts\frac{3}{4}})_2\, B^0_{n + 2,\, n}.
\end{equation}
\smallbreak\item[(ii)]
$\SQ^{\delta-n,\, 1,\, 5}_\lu \cong \SQ^{\delta'-n,\, 1,\, 5}_\lup$ \iff\ they induce simultaneous vanishing of 
\begin{equation} \label{L5p1}
   (\delta - N_5 + {\ts\frac{1}{4}})\, B^1_{n + 3/2,\, n}, \quad
   (\delta - N_5 + {\ts\frac{1}{4}})_2\, B^0_{n + 2,\, n + 1/2}, \quad
   (\delta - N_5 + {\ts\frac{1}{4}})_2\, B^1_{n + 2,\, n}.
\end{equation}
\smallbreak\item[(iii)]
$\SQ^{\delta-n,\, 0,\, 5}_\lu \cong \SQ^{\delta'-n,\, 1,\, 5}_\lup$ \iff\ the first function in~(\ref{L5p0}) at $(\gd)$ and the first function in~(\ref{L5p1}) at $(\gdp)$ are either both zero or both non-zero, and similarly for the second functions in the two displays, and again for the third functions.
\end{enumerate}
\end{prop}

Here $l$ is odd, so the symmetry implied by Corollary~\ref{line dual} appears in the fact that $(N_5, \delta) \mapsto (-N_5, \delta)$ exchanges~(\ref{L5p0}) and~(\ref{L5p1}) (however, in contrast with the condition in Part~(iii) of the proposition, the exchange reorders the functions).  As before, vanishing of the Pochhammer coefficients of the functions $B^\p_{m+r,\, m}$ indicates a splitting due to~(\ref{Dsplitting}).

\subsection{Length $l = 6$} \label{L6}

We have seen that in lengths $l \le 5$, almost all non-resonant \sq s $\SQ^{\delta-n,\, \p,\, l}_\lu$ with a given $n$ are equivalent.  In length~$6$ this is no longer true: for each choice of \cs\ and parity there is a rational invariant whose level curves form a pencil of conics in $(\gd)$-space.  In order to state the results it is helpful to make the following definition.

\meno {\bf Definition.}
Two non-resonant \sq s $\SQ^{\delta-n,\, \p,\, l}_\lu$ and $\SQ^{\delta'-n,\, \p',\, l}_\lup$ are said to satisfy the {\em simultaneous vanishing condition\/} (SVC) if for all pairs $(i,j)$ of elements of $\oh \bN$ such that $i \le \oh (l-1)$ and $\frac{3}{2} \le i - j \le \frac{5}{2}$, one of the following holds:
\begin{enumerate}
\item[(i)]
In the case $\p = \p' = 0$, the \sq s induce simultaneous vanishing of
\begin{equation} \label{SVC0}
   \bigl(\delta - n - \lc j \rc \bigr)_{\lc i \rc - \lc j \rc}\, B^{2\{ j \}}_{n+i,\, n+j}.
\end{equation}
\smallbreak \item[(ii)]
In the case $\p = \p' = 1$, the \sq s induce simultaneous vanishing of
\begin{equation} \label{SVC1}
   \bigl(\delta - n - \oh - \lf j \rf \bigr)_{\lf i \rf - \lf j \rf}\, B^{1 -2\{ j \}}_{n+i,\, n+j}.
\end{equation}
\smallbreak \item[(iii)]
In the case $\p = 0$ and $\p' = 1$, (\ref{SVC0}) at $(\gd)$ and~(\ref{SVC1}) at $(\gdp)$ are either both zero or both non-zero.
\end{enumerate}

\medbreak
Observe that the SVC is vacuous in lengths $l \le 3$, and corresponds to the conditions of Propositions~\ref{L4 prop} and~\ref{L5 prop} in lengths~$4$ and $5$, respectively.  Therefore we have the following unified statement in lengths $l \le 5$.

\begin{prop} \label{Lle5 prop}
For $l \le 5$ and $n$ non-resonant with respect to $l$, the SVC is necessary and sufficient for the equivalence of\/ $\SQ^{\delta-n,\, \p,\, l}_\lu$ and\/ $\SQ^{\delta'-n,\, \p',\, l}_\lup$.
\end{prop}

Now consider length~ $6$.  Here the set of resonant values of $n$ is $\{-2, -\frac{3}{2}, -1, -\oh, 0 \}$, and so that of $N_6$ is $\{0, \pm \oh, \pm 1 \}$.  The SVC involves the functions $B^\p_{n+3/2,\, n}$, $B^{\p+1}_{n+2,\, n+1/2}$, $B^\p_{n+5/2,\, n+1}$, $B^\p_{n+2,\, n}$, $B^{\p+1}_{n+5/2,\, n+1/2}$, and $B^\p_{n+5/2,\, n}$, which we now give in terms of $N_6$.  Recall that $B^0_{m+3/2,\, m}$ is $\gamma^{1/2}$ for all $m$.  For $\p = 0$, the other functions are
\begin{align*}
   B^1_{n+2,\, n+1/2} &=
   \gamma - \bigl[ 3 N_6 \delta + {\ts\frac{3}{4}} \bigr],\ &
   B^0_{n+2,\, n} &=
   \gamma - \bigl[ (N_6 + \oh) (2\delta + N_6 - \oh) \bigr], \\[6pt]
   B^0_{n+5/2,\, n} &=
   \gamma - \bigl[ N_6 \delta + {\ts\frac{3}{4}} \bigr],\ &
   B^1_{n+5/2,\, n+1/2} &=
   \gamma - \bigl[ (N_6 - \oh) (2\delta + N_6 + \oh) \bigr].
\end{align*}
For $\p = 1$, the other functions are
\begin{align*}
   B^1_{n+3/2,\, n} &=
   \gamma - \bigl[ 3(N_6 - \oh) \delta + {\ts\frac{3}{4}} \bigr],\ &
   B^1_{n+2,\, n} &=
   \gamma - \bigl[ (N_6 - 1) (2\delta + N_6) \bigr], \\[6pt]
   B^1_{n+5/2,\, n+1} &=
   \gamma - \bigl[ 3(N_6 + \oh) \delta + {\ts\frac{3}{4}} \bigr],\ &
   B^0_{n+5/2,\, n+1/2} &=
   \gamma - \bigl[ (N_6 + 1) (2\delta + N_6) \bigr],
\end{align*}
\vspace{-6pt}
\begin{equation*}
   B^1_{n+5/2,\, n} =
   \gamma^{3/2} - \gamma^{1/2} \bigl[ 4N_6 \delta - N_6^2 + 3 \bigr].
\end{equation*}

The invariant which determines equivalence is
\begin{equation*}
   I_n^\p(\gd) = B^\p_{n+5/2,\, n}\, B^{\p+1}_{n+2,\, n+1/2}\, \big/\, B^{\p+1}_{n+5/2,\, n+1/2}\, B^\p_{n+2,\, n}.
\end{equation*}

\begin{thm} \label{L6 thm}
For $n$ non-resonant with respect to $l=6$, the SVC is necessary for the equivalence of\/ $\SQ^{\delta-n,\, \p,\, 6}_\lu$ and\/ $\SQ^{\delta'-n,\, \p',\, 6}_\lup$.  Sufficient conditions are as follows:

\begin{enumerate}
\item[(i)]
Suppose that $\p = \p' = 0$.  If either $\delta - n$ or $\delta' - n$ is in $\{0, 1, 2\}$, or if at least one of the factors $B^0_{n+5/2,\, n}$, $B^1_{n+2,\, n+1/2}$, $B^1_{n+5/2,\, n+1/2}$, and $B^0_{n+2,\, n}$ of $I^0_n$ is zero at either $(\gd)$ or $(\gdp)$, then the SVC is also sufficient for equivalence.  Otherwise the \sq s are equivalent \iff\ they satisfy the SVC and in addition
\begin{equation*}
   I^0_n(\gd) = I^0_n(\gdp).
\end{equation*}

\smallbreak\item[(ii)]
Suppose that $\p = \p' = 1$.  If either $\delta - n$ or $\delta' - n$ is in $\{\oh, \frac{3}{2}\}$, or if at least one of the factors $B^1_{n+5/2,\, n}$, $B^0_{n+2,\, n+1/2}$, $B^0_{n+5/2,\, n+1/2}$, and $B^1_{n+2,\, n}$ of $I^1_n$ is zero at either $(\gd)$ or $(\gdp)$, then the SVC is also sufficient for equivalence.  Otherwise the \sq s are equivalent \iff\ they satisfy the SVC and in addition
\begin{equation*}
   I^1_n(\gd) = I^1_n(\gdp).
\end{equation*}

\smallbreak\item[(iii)]
Suppose that $\p=0$ and $\p'=1$.  If either (\ref{SVC0}) at $(\gd)$ or~(\ref{SVC1}) at $(\gdp)$ is zero for $(i,j)$ equal to at least one of $(\frac{5}{2}, 0)$, $(2, \oh)$, $(\frac{5}{2}, \oh)$, and $(2, 0)$, then the SVC is also sufficient for equivalence.  Otherwise the \sq s are equivalent \iff\ they satisfy the SVC and in addition
\begin{equation*}
   I^0_n(\gd) = I^1_n(\gdp).
\end{equation*}
\end{enumerate}
\end{thm}

Let us discuss the implications of this theorem.  Its main point is that together with the SVC, $I^\p_n$ is a complete invariant for the equivalence class of $\SQ^{\delta-n,\, \p,\, 6}_\lu$.  Counting multiplicities, there are exactly four points in the $(\gd)$-plane where the numerator and denominator of $I^\p_n$ are simultaneously zero, and its level curves form the pencil of conics through these four points.  There is an essentially unique modification $\t\gamma_6$ of the coordinate $\gamma$ such that these curves are all in rectilinear orientation in $(\t\gamma_6, \delta)$-coordinates: $\t\gamma_6 := \gamma - 2 N_6 \delta$.  We obtain
\begin{equation*}
   I^0_n = \frac{(\t\gamma_6 - {\ts\frac{3}{4}})^2 - N_6^2 \delta^2}
   {(\t\gamma_6 + {\ts\frac{1}{4}} - N_6^2)^2 - \delta^2}, \qquad
   I^1_n = \frac{\t\gamma_6^2 - 4N_6^2 \delta^2
   + (N_6^2 - 3) (\t\gamma_6 + 2N_6 \delta)}
   {(\t\gamma_6 - N_6^2)^2 - (2\delta + N_6)^2}.
\end{equation*}

In order to simplify computations we now replace $I^\p_n$ by an equivalent invariant $\t I^\p_n$, constructed as were $\t I_n$, $\t J_n$, and $\t M_n$ in Section~6 of \cite{CL13}.  Recall that $B^\p_{m+r,\, m}$ was defined in~(\ref{Bij}) for $r$ in $\{\frac{3}{2},\, 2,\, \frac{5}{2}\}$.  It is implicitly defined in \cite{Co09a} for all $r$ in $\frac{3}{2} + \oh \bN$.  The following lemma may be checked directly from~(\ref{Bij}), and it is no doubt the beginning of an infinite sequence of factorizations analogous to those given in Corollary~7.12 of \cite{CL13} for $\VR$.

\begin{lemma} \label{res facs}
For $\p = 0$ or $1$,
\begin{equation*} \begin{array}{rclrcl}
   B^\p_{2,\, 0} &=& \gamma^{\p/2} B^{\p+1}_{2,\, 1/2}, \quad & \quad
   B^\p_{1/2,\, -3/2} &=& \gamma^{(1-\p)/2} B^\p_{0,\, -3/2}, \\[6pt]
   B^\p_{5/2,\, 0} &=& \gamma^{\p/2} B^{\p+1}_{5/2,\, 1/2}, \quad & \quad
   B^\p_{1/2,\, -2} &=& \gamma^{\p/2} B^\p_{0,\, -2}.
\end{array} \end{equation*}
\end{lemma}

It follows from this lemma that the difference between the numerator and the denominator of $I^\p_n$ is divisible by $N_6^2 - 1$, so we define
\begin{equation*}
   B^\p_{5410} := \bigl( B^\p_{n+5/2,\, n}\, B^{\p+1}_{n+2,\, n+1/2} -
   B^{\p+1}_{n+5/2,\, n+1/2}\, B^\p_{n+2,\, n} \bigr)\, \big/\, \bigl( N_6^2 -1 \bigr),
\end{equation*}
giving
\begin{equation*}
   B^0_{5410} = 2\t\gamma_6 - \delta^2 - N_6^2 - \oh, \qquad
   B^1_{5410} = 3\t\gamma_6 - 4\delta^2 + 2N_6 \delta - N_6^2.
\end{equation*}
The invariants equivalent to $I^0_n$ and $I^1_n$ are
\begin{eqnarray*}
   \t I^0_n &:=& \frac{B^0_{n+5/2,\, n}\, B^1_{n+2,\, n+1/2}}{B^0_{5410}}\ =\
  \frac{(\t\gamma_6 - {\ts\frac{3}{4}})^2 - N_6^2 \delta^2}{2\t\gamma_6 - \delta^2 - N_6^2 - \oh}, \\[6pt]
   \t I^1_n &:=& \frac{B^0_{n+5/2,\, n+1/2}\, B^1_{n+2,\, n}}{B^1_{5410}}\ =\
   \frac{(\t\gamma_6 - N_6^2)^2 - (2\delta + N_6)^2}{3\t\gamma_6 - 4\delta^2 + 2N_6 \delta - N_6^2}.
\end{eqnarray*}
For each $\p$, $I^\p_n$ and $\t I^\p_n$ are related by a linear fractional transformation, so their level curves comprise the same pencil of conics.

For $\p = 0$ the pencil is defined by the four points $\bigl(N_6^2 + \epsilon_1 N_6 +1,\, \epsilon_2 (N_6 + \epsilon_1) \bigr)$ in $(\t\gamma_6, \delta)$-space, where $\epsilon_\bullet = \pm 1$.  These points form an isosceles trapezoid for all $N_6$.

For $\p = 1$ the pencil is defined by the two points $\bigl( N_6^2 + 3\epsilon_1 N_6 +3,\, N_6 + \frac{3}{2} \epsilon_1 \bigr)$, where $\ep_1 = \pm 1$, and the double point $\bigl( N_6^2, -\oh N_6 \bigr)$, again in $(\t\gamma_6, \delta)$-space (we have not worked out the tangency condition associated to the double point).  The fact that these two pencils are rational in $N_6$ may have an underlying explanation akin to that given for the invariant $I_n$ in Section~6.2 of \cite{CL13}.

\subsection{Results of \cite{BBB13}}  \label{BBB Results}
This article gives the \ec es of the \dog\ modules $\D^k_\lu$ in all cases, resonant as well as non-resonant.  We now use our results to recover these \ec es in the non-resonant cases.  As noted in Section~\ref{Sbqnts}, $\D^k_\lu$ is the \sq\ $\SQ^{k,\, 2k,\, 2k+1}_\lu$ of length $2k+1$, and $\delta$ is a complete invariant for its \cs.  Since $n = \delta - k$ and $k$ is fixed, one checks that the resonant values of $\delta$ are $\oh,\, 1,\, \frac{3}{2},\, \ldots,\, k$.

In lengths~$1$, $2$, and~$3$, we saw in Section~\ref{Lle5} that non-resonant \sq s split completely.  Hence for $k=0$, $\oh$, or~$1$, $\delta$ is a also complete invariant for the \ec\ of $\D^k_\lu$.

In length~$4$, Proposition~\ref{L4 prop}(ii) shows that $\D^{3/2}_\lu$ is split \iff~(\ref{L4p1}) is zero.  Here $N_4 = \delta - 1$, so $(\delta - N_4) B^1_{n+3/2,\, n}$ reduces to $12 \lambda (\mu - \oh)$.  Thus for any fixed non-resonant $\delta$ there are two \ec es: the split class, where $\lambda$ is either~$0$ or $\oh - \delta$, and the non-split class, containing all the other modules.  Observe that the split class consists of the modules $\D^{3/2}_{0,\, \delta}$ and $\D^{3/2}_{-\delta + 1/2,\, 1/2}$ occurring in Lemma~\ref{DBol}.

In length~$5$, Proposition~\ref{L5 prop}(i) shows that $\D^2_\lu$ and $\D^2_\lup$ are equivalent \iff\ they induce simultaneous vanishing of the three functions~(\ref{L5p0}).  Here $N_5 = \delta - \frac{5}{4}$, so the Pochhammer symbols are all non-zero (as always on \dog\ modules) and one finds that $B^1_{n+2,\, n+1/2}$ and $B^2_{n+2,\, n}$ are both $12 \lambda (\mu - \oh)$.  Therefore the condition reduces to simultaneous vanishing of $\lambda (\mu - \oh)$ and $\lambda + \mu - \oh$, so for each non-resonant $\delta$ there are three equivalence classes: $\lambda = \oh (\oh - \delta)$ gives the self-adjoint class; $\lambda = 0$ and $\lambda = \oh - \delta$ form the class in which the order~0 operators split off by Lemma~\ref{DBol}; and all the other values of $\lambda$ comprise the minimally split class.

We remark that $\lambda (\mu - \oh)$ and $\lambda + \mu - \oh$ both vanish only on the self-adjoint module $\D^2_{0,\, 1/2}$.  Although $\delta = \oh$ is resonant, it follows from Lemma~\ref{DBol} that this module is completely split.

In length~$6$ we apply Theorem~\ref{L6 thm}(ii) to $\D^{5/2}_\lu$.  Here $N_6 = n+1 = \delta - \frac{3}{2}$, and the factors of the invariant $I^1_{\delta - 5/2}$ reduce to
\begin{equation*} \begin{array}{ll}
   B^1_{n+5/2,\, n} = 12 \gamma^{1/2} \lambda (\mu - \oh), \qquad &
   B^0_{n+2,\, n+1/2} = \gamma^{1/2}, \\[6pt]
   B^0_{n+5/2,\, n+1/2} = 12\lambda (\mu - \oh), \qquad &
   B^1_{n+2,\, n} = \gamma - 3(\delta - \oh) (\delta - {\ts\frac{5}{2}}).
\end{array} \end{equation*}
By the SVC, the conjugate pair of $\lambda$-values~$0$ and $\oh - \delta$ forms one equivalence class, and away from this class we have
\begin{equation*}
   I^1_{\delta-5/2} = \gamma\, \big/\, \bigl[ \gamma - 3(\delta - \oh) (\delta - {\ts\frac{5}{2}}) \bigr].
\end{equation*}
Since $\delta$ is fixed, the $\gamma$-monotonicity of this function shows that conjugation is the only equivalence: the non-resonant \ec es are the conjugate pairs and the self-adjoint singleton $\lambda = \oh (\oh - \delta)$.

The length~$6$ result immediately extends to all higher lengths, as an equivalence between $\D^k_\lu$ and $\D^k_\lup$ must restrict to an equivalence between $\D^{k'}_\lu$ and $\D^{k'}_\lup$ for all $k' < k$.  Thus for $k \ge \frac{5}{2}$, non-resonant modules $\D^k_\lu$ and $\D^k_\lup$ are equivalent \iff\ they are equal or conjugate, \ie\ $(\gd) = (\gdp)$.  This matches the non-resonant results of \cite{BBB13}.

\section{Lengths $l \ge 7$}  \label{Lge7}

In this section we state partial results in higher lengths.  As in Section~\ref{Lle6}, the proofs are postponed to Section~\ref{Proofs}.  We begin with a general result.

\begin{prop} \label{SVC prop}
In all lengths~$l$, the SVC is necessary for\/ $\SQ^{\delta-n,\,\p,\, l}_\lu \cong \SQ^{\delta'-n,\, \p',\, l}_\lup$.
\end{prop}

Henceforth we shall treat only even equivalences, so we shall be concerned only with the sets $\EC^{\p,l}_n(\gd)$.  In odd lengths it suffices to compute $\EC^{0,l}_n(\gd)$, because, writing $\bigl(\begin{smallmatrix} 1 & 0 \\ 0 & -1 \end{smallmatrix}\bigr)$ for the transformation $(\gd) \mapsto (\gamma, -\delta)$, Corollary~\ref{line dual} gives
\begin{equation} \label{odd length}
   \EC^{1,l}_n(\gd) = \bigl(\begin{smallmatrix} 1 & 0 \\ 0 & -1 \end{smallmatrix}\bigr) \EC^{0,l}_{1-n-l/2}(\gamma, -\delta).
\end{equation}

\subsection{Length $l = 7$} \label{L7}

Since $\SQ^{\delta-n,\, \p,\, 7}_\lu$ contains the two \sq s $\SQ^{\delta-n,\, \p,\, 6}_\lu$ and $\SQ^{\delta-n-1/2,\, \p+1,\, 6}_\lu$ of length~$6$, it is clear that $I^\p_n$ and $I^{\p+1}_{n+1/2}$ are both invariants.  There is a third invariant:
\begin{equation*}
   J^\p_n(\gd)\ :=\ B^{\p+1}_{n+3,\, n+1/2}\, B^\p_{n+5/2,\, n}\, \big/\,
   B^{\p+1}_{n+3,\, n+3/2}\, B^{\p+1}_{n+5/2,\, n+1/2}\, B^\p_{n+3/2,\, n}.
\end{equation*}

Observe that the numerator and denominator of $J^\p_n$ both have $\gamma^{1/2}$ as a factor, and cancelling it gives a ratio of conics.  However, the value of this simplified form of $J^\p_n$ has no meaning at $\gamma = 0$, and $J^\p_n$ should be regarded as undefined there.

In contrast with this, $I^\p_n$, $I^{\p+1}_{n+1/2}$, and $J^\p_n$ share various factors $B^\bullet_{\bullet, \bullet}$, and when some of these factors are zero it can happen that some of the invariants themselves are undefined but some of their ratios are defined after cancelling these factors, and moreover, the values of said ratios are meaningful.  For this reason, if $\I$ is a ratio of products of invariants, we define $\Simp(\I)$ to be the ratio of products of factors $B^\bullet_{\bullet, \bullet}$ obtained by cancelling as many such factors as possible.

By~(\ref{odd length}), $\EC^{1,7}_n(\gd)$ is the reflection of $\EC^{0,7}_{-n-5/2}(\gamma, -\delta)$, so we only describe the classes $\EC^{0, 7}_n$.

\begin{thm} \label{L7 thm}
For $n$ non-resonant with respect to $l=7$,\/ $\SQ^{\delta-n,\, 0,\, 7}_\lu \cong \SQ^{\delta'-n,\, 0,\, 7}_\lup$ \iff\ they satisfy the SVC and at least one of the following conditions:

\begin{enumerate}
\item[(i)]
At least one of $\delta-n$ and $\delta'-n$ is~$1$ or~$2$.

\smallbreak\item[(ii)]
At least one of $\delta-n$ and $\delta'-n$ is~$0$, and\/ $\SQ^{\delta-n-1/2,\, 1,\, 6}_\lu \cong \SQ^{\delta'-n-1/2,\, 1,\, 6}_\lup$.

\smallbreak\item[(iii)]
For all ratios $\I$ of products of $I^0_n$, $I^1_{n+1/2}$, and $J^0_n$ such that none of the factors $B^\bullet_{\bullet, \bullet}(\gd)$ of\/ $\Simp(\I)$ is zero,
\begin{equation} \label{Simp}
   \Simp\bigl(\I(\gd)\bigr) = \Simp\bigl(\I(\gdp)\bigr).
\end{equation}
\end{enumerate}
\end{thm}

Parts~(i) and~(ii) of this theorem reflect~(\ref{Dsplitting}).  Let us explain the somewhat mysterious condition in Part~(iii).  Note that $I^0_n$, $I^1_{n+1/2}$, and $J^0_n$ have a total of nine distinct factors $B^{2j}_{n+i,\, n+j}$: those with $i, j \in \oh\bN$, $i \le 3$, and $\frac{3}{2} \le i-j \le \frac{5}{2}$.  The recurring factors are as follows: $B^0_{n+5/2,\, n}$ is shared by $I^0_n$ and $J^0_n$, $B^1_{n+3,\, n+1/2}$ is shared by $I^1_{n+1/2}$ and $J^0_n$, and $B^1_{n+5/2,\, n+1/2}$ is shared by all three.

In general there are infinitely many allowed ratios $\I$ in Condition~(iii), but in fact it is never necessary to consider more than three of them.  As the number of factors vanishing at $(\gd)$ and $(\gdp)$ increases, the number of ``basic'' ratios $\I$ which must be considered decreases.  The simplified ratios which arise are
\begin{eqnarray*}
   \Simp\Bigl(\frac{J^0_n}{I^0_n}\Bigr)
   &:=& \frac{B^1_{n+3,\, n+1/2}\, B^0_{n+2,\, n}}
   {B^1_{n+3,\, n+3/2}\, B^1_{n+2,\, n+1/2}\, B^0_{n+3/2,\, n}},
   \\[6pt]
   \Simp\Bigl(\frac{J^0_n}{I^1_{n+1/2}}\Bigr)
   &:=& \frac{B^0_{n+5/2,\, n}\, B^0_{n+3,\, n+1}}
   {B^1_{n+3,\, n+3/2}\, B^0_{n+5/2,\, n+1}\, B^0_{n+3/2,\, n}},
   \\[6pt]
   \Simp\Bigl(\frac{J^0_n}{I^0_n I^1_{n+1/2}}\Bigr)
   &:=& \frac{B^0_{n+3,\, n+1}\, B^1_{n+5/2,\, n+1/2}\, B^0_{n+2,\, n}}
   {B^1_{n+3,\, n+3/2}\, B^0_{n+5/2,\, n+1}\, B^1_{n+2,\, n+1/2}\, B^0_{n+3/2,\, n}}.
\end{eqnarray*}

In the following table, some of the various possible sets of vanishing factors $B^\bullet_{\bullet, \bullet}$ are listed on the left.  In each case, Condition~(iii) holds \iff~(\ref{Simp}) holds for each invariant on the right.

\renewcommand{\arraystretch}{1.5}
\bigbreak \begin{center} \begin{tabular}{|l|l|} \hline
{\bf Vanishing} $B^\bullet_{\bullet, \bullet}$ & {\bf Basic Invariants} \\ \hline
None &
$I^0_n$, $I^1_{n+1/2}$, and $J^0_n$ \\ \hline
$B^0_{n+3/2,\, n}$ and/or $B^1_{n+3,\, n+3/2}$ &
$I^0_n$ and $I^1_{n+1/2}$ \\ \hline
$B^0_{n+5/2,\, n+1}$ and/or $B^0_{n+3,\, n+1}$ $\qquad$ &
$I^0_n$ and $J^0_n$ \\ \hline
$B^0_{n+2,\, n}$ and/or $B^1_{n+2,\, n+1/2}$  &
$I^1_{n+1/2}$ and $J^0_n$ \\ \hline
$B^1_{n+5/2,\, n+1/2}$ &
$\Simp(J^0_n / I^0_n)$ and $\Simp(J^0_n / I^1_{n+1/2})$ $\qquad$ \\ \hline
$B^0_{n+5/2,\, n}$ &
$\Simp(J^0_n / I^0_n)$ and $I^1_{n+1/2}$ \\ \hline
$B^1_{n+3,\, n+1/2}$ &
$\Simp(J^0_n / I^1_{n+1/2})$ and $I^0_n$ \\ \hline
$B^0_{n+5/2,\, n}$ and $B^1_{n+3,\, n+1/2}$ &
$\Simp(J^0_n / I^0_n I^1_{n+1/2})$ \\ \hline
$B^0_{n+5/2,\, n}$  and $B^1_{n+5/2,\, n+1/2}$ &
$\Simp(J^0_n / I^0_n)$ \\ \hline
$B^1_{n+3,\, n+1/2}$ and $B^1_{n+5/2,\, n+1/2}$ &
$\Simp(J^0_n / I^1_{n+1/2})$ \\ \hline
\end{tabular} \end{center} \bigbreak
\renewcommand{\arraystretch}{1}

It is also possible for any of $I^0_n$, $I^1_{n+1/2}$, and $J^0_n$ to be the sole basic invariant, or, if enough factors vanish, for Condition~(iii) to be vacuous.

We have given formulas for $I^0_n$ and $I^1_n$ in terms of $N_6$, from which formulas for $I^0_n$ and $I^1_{n+1/2}$ in terms of $N_7 = N_6 + \frac{1}{4}$ follow immediately.  The formula for $J^0_n$ may easily be derived from~(\ref{Bij}).  We give it in terms $N_6$ because there is no clear advantage in using $N_7$: for $\gamma \not= 0$,
\begin{equation*}
   J^0_n\ =\ \frac{ \bigl[ \gamma - 2(2N_6 + 1) \delta + (N_6 + \oh)^2 - 3 \bigr]\,
   \bigl[ \gamma - N_6 \delta - {\ts\frac{3}{4}} \bigr] }
   { \bigl[ \gamma - (N_6 - \oh) (2\delta + N_6 + \oh) \bigr]\,
   \bigl[ \gamma - 3(N_6 + 1) \delta - {\ts\frac{3}{4}} \bigr] }.
\end{equation*}

Just as for $I^p_n$, the level curves of this function form a pencil of conics.  As in length~$6$, there is a choice of coordinate $\t\gamma_7$ such that all conics in the pencil are rectilinear in $(\t\gamma_7, \delta)$-coordinates: $\t\gamma_7 := \gamma - (\frac{5}{2} N_6 + 1) \delta$.  We remark that there are similarities between the coordinates $\t\gamma_6$ and $\t\gamma_7$ defined in this article and the coordinates $\t\gamma_5$ and $\t\gamma_6$ defined in Section~6 of \cite{CL13}, but we do not know an underlying explanation of the pattern.

Also as in length~$6$, there are special values of $N_6$ where $J^0_n$ reduces to~$1$, namely $N_6 = 1$ and $-\frac{3}{2}$, \ie\ $N_7 = \pm\frac{5}{4}$.  We do know an explanation of this phenomenon: in Section~\ref{Resonance} we will describe an extension of~(\ref{Bij}) defining $B^\p_{m+r,\, m}$ for all $r \in \frac{3}{2} + \oh \bN$, and the {\em cup equation\/} (see Section~4.2 of \cite{Co09a}) implies that $B^0_{n+3,\, n}$ is equal to $B^1_{n+3,\, n+3/2} B^0_{n+3/2,\, n}$.  Coupling this with the extension of Lemma~\ref{res facs} to $B^0_{n+3,\, n}$ gives the reduction.  Following our construction of $\t I^\p_n$, this reduction may be used to replace $J^0_n$ by a simpler invariant $\t J^0_n$ with the same level curves:
\begin{eqnarray*}
   \t J^0_n &:=& \frac
   { 2(N_6 - 1) (N_6 + {\ts\frac{3}{2}})\,
   B^1_{n+3,\, n+3/2}\, B^1_{n+5/2,\, n+1/2}\, B^0_{n+3/2,\, n} }
   { B^1_{n+3,\, n+1/2}\, B^0_{n+5/2,\, n}\, -\,
   B^1_{n+3,\, n+3/2}\, B^1_{n+5/2,\, n+1/2}\, B^0_{n+3/2,\, n} } \\[6pt]
   &=& \frac { \bigl[ \gamma - (N_6 - \oh) (2\delta + N_6 + \oh) \bigr]\,
   \bigl[ \gamma - 3(N_6 + 1) \delta - {\ts\frac{3}{4}} \bigr] }
   { \gamma - \delta^2 - (2N_6 + 1) \delta - {\ts\frac{3}{4}} }.
\end{eqnarray*}

\meno {\bf Remark.}
At the beginning of Section~\ref{Lge7} we described our results in lengths exceeding~$6$ as partial.  In length~$7$ this is because Theorem~\ref{L7 thm} is clearly inconclusive, even though its conditions of are necessary and sufficient for equivalence: it states that generically the \sq s are equivalent \iff\ $(\gd)$ and $(\gdp)$ lie on the same conic in each of three pencils, and generically three conics intersect only in a single point.  Therefore one expects that for almost all choices of $(n, \lu)$, the even \ec\ $\EC^{0,7}_n(\gd)$ consists only of the single point $(\gd)$, \ie\ $\SQ^{k, 0, 7}_\lu$ is equivalent only to its conjugate, while possibly for finitely many special values of $n$ (which could turn out to be resonant, invalidating them) it consists of two, three, or four points.

The obvious attempt to resolve this consists essentially in solving for $(\gd)$ in terms of $\t I^0_n$, $\t I^1_{n+1/2}$, and $\t J^0_n$.  It is possible to reduce this to a Gr\"obner basis problem amenable to software, but we made only preliminary explorations in this direction.

\subsection{Lengths $l \ge 8$}  \label{Lge8}

The following result shows that the rational invariants occurring in length~$7$ completely govern even equivalence in higher lengths, except in the case of certain exceptional \cs\ where there may be one additional invariant.

\begin{prop} \label{Lge8 prop}
Fix\/ $l \ge 8$ and $n$ non-resonant with respect to $l$.  Suppose first that neither of\/ $- n - \frac{1}{4}(7 \pm \sqrt{33})$ is an element of\/ $\oh \bN$ less than $\oh l - 4$, or in other words that the composition series of\/ $\SQ^{\delta-n,\, \p,\, l}_\lu$ does not contain either of the following two pairs  of exceptional \tdm s in either parity:
\begin{equation} \label{jump 4}
   \bigl\{ \F_{(-7+\sqrt{33})/4},\, \F_{(9+\sqrt{33})/4} \bigr\}, \qquad
   \bigl\{ \F_{(-7-\sqrt{33})/4},\, \F_{(9-\sqrt{33})/4} \bigr\}.
\end{equation}
Then\/ $\SQ^{\delta-n,\, \p,\, l}_\lu \cong \SQ^{\delta'-n,\, \p,\, l}_\lup$ \iff\ all of their \sq s of length~$7$ are equivalent pairwise, or in other words \iff\ for all $i$ in $\oh \bN$ less than $\oh l - 3$,
\begin{equation*}
   \SQ^{\delta - n - i,\, \p + 2i,\, 7}_\lu\ \cong\
   \SQ^{\delta' - n - i,\, \p + 2i,\, 7}_\lup.
\end{equation*}

If on the other hand one of\/ $- n - \frac{1}{4}(7 \pm \sqrt{33})$ is an element $i$ of\/ $\oh \bN$ less than $\oh l - 4$, then\/ $\SQ^{\delta-n,\, \p,\, l}_\lu \cong \SQ^{\delta'-n,\, \p,\, l}_\lup$ \iff\ all of their \sq s of length~$7$ are equivalent pairwise and in addition,
\begin{equation*}
   \SQ^{\delta - n - i,\, \p + 2i,\, 9}_\lu\ \cong\
   \SQ^{\delta' - n - i,\, \p + 2i,\, 9}_\lup.
\end{equation*}
\end{prop}

The point of this proposition is that when the pairs~(\ref{jump 4}) are not present, the rational invariants
\begin{equation} \label{arb length}
   I_n^\p,\, I_{n+1/2}^{\p+1},\, \ldots, I_{n+(l-6)/2}^{\p+l}, \quad
   J_n^\p,\, J_{n+1/2}^{\p+1},\, \ldots, J_{n+(l-7)/2}^{\p+l-1}
\end{equation}
are essentially complete, in the same sense that $I_n^0$, $I_{n+1/2}^1$, and $J_n^0$ are complete in Theorem~\ref{L7 thm}(iii): one needs also the SVC and equality of the simplified ratios of all products of the invariants~(\ref{arb length}).

The pairs~(\ref{jump 4}) are exceptional because of the 1-cocycles discovered in Theorem~7.5(f) of \cite{Co09a}.  When one of them is present one must add a single additional invariant to~(\ref{arb length}).  In terms of the function $B^\p_{m+4,\, m}$ we will discuss in Section~\ref{Resonance}, this invariant is
\begin{equation*}
   K^\p_m (\gd)\ :=\ B^\p_{m+4,\, m}\, \big/\, B^\p_{m+4,\, m+2}\, B^\p_{m+2,\, m}, \quad
   m = {\ts\frac{1}{4}} (-7 \pm \sqrt{33}).
\end{equation*}

In fact, we expect that these exceptional \cs\ are an artifact of the incompleteness of our results in lengths exceeding~$6$.  Recall that in length~$7$ there are three rational invariants, each of whose level curves is a pencil of conics, so we predicted that conjugation is the only even equivalence for almost all $(n, \lu)$.  In length~$8$ there are five such invariants, so we conjecture that conjugation is the only even equivalence in all cases, \ie\ $\EC^{\p, 8}_n(\gd)$ is always simply $\{(\gd)\}$.  We also conjecture that the only odd equivalences are given by Lemma~\ref{Bol}.

A relatively simple computation using software shows that provided the \sq s do not split as in~(\ref{Dsplitting}), this is true for even equivalences in lengths $l \ge 15$.  We emphasize that $l = 15$ is probably not special; a more detailed analysis should give a lower length.

\begin{prop} \label{Lge15 prop}
Fix $n$ non-resonant with respect to $l = 15$.  For $\p = 0$, assume that neither $(\delta - n)_6$ nor $(\delta' - n)_6$ is zero, and for $\p = 1$, assume that neither $(\delta - n - \oh)_6$ nor $(\delta' - n - \oh)_6$ is zero.  Then\/ $\SQ^{\delta-n,\, \p,\, 15}_\lu \cong \SQ^{\delta'-n,\, \p,\, 15}_\lup$ \iff\ they are either equal or conjugate.
\end{prop}

\section{Lacunary \sq s}  \label{Lacunarity}

By Proposition~7.1 of \cite{Co09a}, if $n = \delta - k$ is neither~$0$ nor~$\oh$ then $\Psi^{k,\, \p}_\lu$ contains a unique {\em lacunary\/} $\K$-submodule $\Psi^{k,\, \p,\, \lac}_\lu$, which contains $\Psi^{k-3/2,\, \p+1}_\lu$ and has \cs
\begin{equation*}
   \bigl\{ \F^{\p\Pi}_n,\, \F^{(\p+1)\Pi}_{n+3/2},\, \F^{p\Pi}_{n+2},\, \F^{(p+1)\Pi}_{n+5/2},\, \ldots \bigr\}.
\end{equation*}
It is constructed by taking the inverse image of $\F^{\p\Pi}_n$ under the map $\Symb_{3/2}$ given in~(25) of \cite{Co09a}, and may be thought of as the space of \psdog s from $\F_\lambda$ to $\F_\mu$ of order~$\le (k, \p)$ without order $k-\oh$ or $k-1$ terms.  The point is that there is a $\K$-invariant way to specify the order $k-\oh$ and $k-1$ terms of an order $k$ operator: in general this is not true for terms of order~$\le k-\frac{3}{2}$.  It is of course also possible to remove only the order $k-\oh$ or only the $k-1$ terms: the resulting modules are $\Psi^{k,\, \p,\, \lac}_\lu + \Psi^{k-1,\, \p}_\lu$ and $\Psi^{k,\, \p,\, \lac}_\lu + \Psi^{k-1/2,\, \p+1,\, \lac}_\lu$, respectively.

The equivalence question may be posed for \sq s of lacunary modules.  We will only consider the ``maximally lacunary \sq''
\begin{equation*}
   \SQ^{k,\, \p,\, l,\, \lac}_\lu\ :=\ \Psi^{k,\, \p,\, \lac}_\lu\, \big/\,
   \bigl( \Psi^{k-(l+1)/2,\, \p+l+1,\, \lac}_\lu\, +\, \Psi^{k-(l-2)/2,\, \p+l,\, \lac}_\lu \bigr).
\end{equation*}
For $l \ge 3$, this module is of length~$l$ with \cs\
\begin{equation*}
   \bigl\{ \F^{\p\Pi}_n,\, \F^{(\p+1)\Pi}_{n+3/2},\, \F^{\p\Pi}_{n+2},\, \ldots,\,
   \F^{(\p+l+1)\Pi}_{n+(l-1)/2},\, \F^{(\p+l)\Pi}_{n+l/2},\, \F^{(\p+l+1)\Pi}_{n+(l+3)/2} \bigr\}.
\end{equation*}

We shall say that $\SQ^{\delta - n,\, \p,\, l,\, \lac}_\lu$ and $\SQ^{\delta'-n,\, \p,\, l,\, \lac}_\lup$ satisfy the {\em lacunary SVC\/} for $\p = 0$ (respectively, $\p = 1$) if they induce simultaneous vanishing of the functions~(\ref{SVC0}) (respectively, (\ref{SVC1})) for all pairs $(i,j)$ of elements of $\oh \bN$ such that $j=0$ or $j \ge \frac{3}{2}$, $i \le \oh l$ or $i = \oh (l+3)$, and $\frac{3}{2} \le i-j \le \frac{5}{2}$.

At $l=4$ the lacunary SVC involves four pairs $(i,j)$, namely,
\begin{equation*}
   ({\ts\frac{7}{2}}, {\ts\frac{3}{2}}), \qquad ({\ts\frac{3}{2}}, 0), \qquad 
   ({\ts\frac{7}{2}}, 2), \qquad (2, 0).
\end{equation*}
Here there is a single rational invariant:
\begin{equation*}
   M^\p_n(\gd)\ :=\ B^{\p+1}_{n+7/2,\, n+3/2}\, B^\p_{n+3/2,\, n}\, \big/\, B^\p_{n+7/2,\, n+2}\, B^\p_{n+2,\, n}.
\end{equation*}

\begin{prop} \label{lac 4 prop}
For $l \ge 3$ and $n$ non-resonant with respect to $l$, the lacunary SVC is necessary for\/ $\SQ^{\delta - n,\, \p,\, l,\, \lac}_\lu \cong \SQ^{\delta'-n,\, \p,\, l,\, \lac}_\lup$.  For $l = 3$ it is also sufficient.

For $l = 4$, if any of the four functions involved in the SVC is zero then the SVC is still sufficient for equivalence.  Otherwise the \sq s are equivalent \iff\
\begin{equation*}
   M^\p_n(\gd) = M^\p_n(\gdp).
\end{equation*}
\end{prop}

The invariant $M^\p_n$ is best expressed in terms of $N_8 = n + \frac{3}{2}$.  For $\p=0$ and $\gamma \not= 0$, it simplifies to a ratio of linear functions:
\begin{equation*}
   M^0_n := \bigl[ \gamma - 2N_8 \delta - N_8^2 - N_8 \bigr]\, \big/\,
   \big[( \gamma - 2N_8 \delta - N_8^2 + N_8 \bigr].
\end{equation*}
Clearly we may replace $M^0_n$ by the equivalent invariant $\gamma - 2N_8 \delta$.

At $\p = 1$ we find
\begin{equation*}
   M^1_n\ :=\ \frac{\bigl[ \gamma - (2N_8 + 3) \delta - (N_8 + \oh) (N_8 + {\ts\frac{3}{2}}) \bigr]\,
   \bigl[ \gamma - 3(N_8 - 1) \delta - {\ts\frac{3}{4}} \bigr]}
   {\bigl[ \gamma - (2N_8 - 3) \delta - (N_8 - \oh) (N_8 - {\ts\frac{3}{2}}) \bigr]\,
   \bigl[ \gamma - 3(N_8 + 1) \delta - {\ts\frac{3}{4}} \bigr]}.
\end{equation*}
It turns out that setting $\t\gamma_8 := \gamma - \frac{5}{2} N_8 \delta$ gives a rectilinear pencil of conic level curves, following the (as yet unexplained) pattern set by $\t\gamma_6$ and $\t\gamma_7$.  Lemma~\ref{res facs} shows that $M^1_n = 1$ at $N_8=0$, so we may replace it with the simpler invariant
\begin{eqnarray*}
   \t M^1_n &:=& \frac{ -2N_8
   \bigl[ B^0_{n+7/2,\, n+3/2}\, B^1_{n+3/2,\, n} + B^1_{n+7/2,\, n+2}\, B^1_{n+2,\, n} \bigr]}
   {B^0_{n+7/2,\, n+3/2}\, B^1_{n+3/2,\, n} - B^1_{n+7/2,\, n+2}\, B^1_{n+2,\, n}} \\[6pt]
   &=& \frac{ 4 \t\gamma_8^2 - (N_8^2 + 36) \delta^2 - 2(2N_8^2 + 3) \t\gamma_8
   + 2N_8 (N_8^2 - 12) \delta + 3(N_8^2 + {\ts\frac{3}{4}}) }
   { 4\t\gamma_8 - 6\delta^2 + 4N_8 \delta - 3 }.
\end{eqnarray*}

We did not investigate the lacunary cases further, except to note that at both $l=5$ and $l=6$ there are two rational invariants, with the possible exception of the length~$5$ case $n = \frac{1}{4} (-7 \pm \sqrt{33})$, where $K^\p_n$ may be a third invariant.

\section{Proofs}  \label{Proofs}

Our proofs are based on \cite{Co09a}.  It is proven in Corollary~6.6 of that paper that if $n$ is non-resonant with respect to $l$, then there exists a unique even $\ds$-equivalence
\begin{equation*}
   \CQ_\lu : \biggl(\th \bigoplus_{i \in \frac{1}{2}\bN}^{2i \le l-1} \F^{(\p+2i)\Pi}_{n+i} \th\biggr)\,
   \to\, \SQ^{\delta-n,\, \p,\, l}_\lu
\end{equation*}
which {\em preserves symbols:\/} it injects $\F^{(\p+2i)\Pi}_{n+i}$ to the image of $\Psi^{\delta-n-i,\, \p+2i}_\lu$, and composing this injection with the natural symbol map defined in Lemma~\ref{Psi fine filt}(ii) gives the identity.  (In \cite{Co09a}, $\ds$ was called the conformal subalgebra of $\K$ and $\CQ_\lu$ was called the conformal quantization, but this was a misnomer: ``conformal'' should be replaced with ``projective''.)  Using $\CQ_\lu$ we can write the action $L_\lu$ of $\K$ on $\SQ^{\delta-n,\, \p,\, l}_\lu$ in an explicitly $\ds$-diagonal manner:

\meno {\bf Definition.}
Let $\pi^\lu$ be the \r\ of $\K$ on $\bigoplus_{2i=0}^{l-1} \F^{(\p+2i)\Pi}_{n+i}$ given by
\begin{equation*} 
   \pi^\lu(X) := \CQ_\lu^{-1} \circ \big(L_\lu(X) \big|_{\SQ^{\delta-n,\, \p,\, l}_\lu}\big) \circ \CQ_\lu.
\end{equation*}
Regard $\pi^\lu$ as an $l \times l$ matrix with entries
\begin{equation*}
   \pi^{\lu,\, \p+2j \mod 2}_{n+i,\, n+j}: \K \to \Hom(\F^{(\p+2j)\Pi}_{n+j}, \F^{(\p+2i)\Pi}_{n+i}).
\end{equation*}

\medbreak
The next lemma is drawn from Lemma~5.3 and Theorem~6.5 of \cite{Co09a}.  It is a consequence of the $\K$-invariance of the refined order filtration and the $\ds$-covariance of $\CQ_\lu$.  The lemma that follows it is the second paragraph of Lemma~5.4 of \cite{Co09a}.  The subsequent corollary, which is obvious from the two lemmas, is part of Theorem~6.5 of that paper.

\begin{lemma}
The matrix entry $\pi^{\lu,\, \p+2j}_{n+i,\, n+j}$ is zero for $i < j$.  For $i = j$ it is the tensor density action $L^{(\p+2i) \Pi}_{n+i}$.  For $i > j$ it is even, $\ds$-covariant, and zero on $\ds$.
\end{lemma}

\begin{lemma}
The space of even $\ds$-covariant maps from $\K$ to\/ $\Hom(\F^{q \Pi}_m, \F^{q' \Pi}_{m'})$ which vanish on $\ds$ is zero unless $m' - m \in \frac{3}{2} + \oh \bN$ and $q' - q \equiv 2(m' - m)$ modulo~$2$.  In this case it is 1-dimensional and spanned by $\beta_{m, m'}$, its unique element such that
\begin{equation*}
   \beta_{m, m'}(X_{\xi x^2}) = 2 \alpha^{m'-m}
   \o D^{2(m' - m) - 3} \circ \ep_{\F_m}^{2(m'-m)}.
\end{equation*}
\end{lemma}

\begin{cor} \label{pi b's}
For $i - j = \oh$ or~$1$, $\pi^{\lu,\, \p+2j}_{n+i,\, n+j} = 0$.  For $i - j \in \frac{3}{2} + \oh \bN$, there are scalars $b^{\p + 2j}_{n+i,\, n+j}(\lu)$ such that
\begin{equation*}
   \pi^{\lu,\, \p+2j}_{n+i,\, n+j} = b^{\p + 2j}_{n+i,\, n+j} (\lu) \beta_{n+j,\, n+i}.
\end{equation*}
\end{cor}

We remark that in cohomological terms, $\beta_{m, m'}$ is an {\em $\ds$-relative 1-cochain.\/}  We will only need the explicit formula for $b^\p_{m+r,\, m}$ at $r = \frac{3}{2}$, $2$, and $\frac{5}{2}$.  Recall the scalars $B^\p_{m+r,\, m}$ from~(\ref{Bij}).  The following result is given for \sq s of \dog s and then more generally for \sq s of \psdog s in Theorems~5.6 and~6.5 of \cite{Co09a}, respectively.

\begin{thm} \label{bCB}
For $r \in \{ \frac{3}{2},\, 2,\, \frac{5}{2} \}$, $b^\p_{m+r,\, m} = C^\p_{m+r,\, m}\, B^\p_{m+r,\, m}$, where
\begin{align*}
   C^0_{m+3/2,\, m} &= -\frac{1}{4 \sqrt{3}}\, \frac{(\delta - m)_2}{m + 1/2}, &
   C^1_{m+3/2,\, m} &= -\frac{1}{12}\, \frac{\delta - m - 1/2}{m + 1/2} \\[6pt]
   C^0_{m+2,\, m} &= \frac{1}{64}\, \frac{(\delta - m)_2}{m(m + 3/2)}, &
   C^1_{m+2,\, m} &= -\frac{1}{16}\, \frac{(\delta - m - 1/2)_2}{m(m + 3/2)} \\[6pt]
   C^0_{m+5/2,\, m} &= \frac{1}{16}\, \frac{(\delta - m)_3}{(m + 2)_3}, &
   C^1_{m+5/2,\, m} &= \frac{1}{16 \sqrt{3}}\, \frac{(\delta - m - 1/2)_2}{(m + 2)_3}
\end{align*}
\end{thm}

Observe that the denominators of the scalars $C^{\p+2j}_{n+i,\, n+j}$ have zeroes precisely in the resonant cases.  This theorem and the following general condition under which two non-resonant \sq s are equivalent comprise the key to our approach.

\begin{prop} \label{nonres eqvs}
For $l$ arbitrary and $n$ non-resonant with respect to $l$, the \sq s\/ $\SQ^{\delta-n,\, \p,\, l}_\lu$ and\/ $\SQ^{\delta'-n,\, \p',\, l}_\lup$ are equivalent (necessarily by an equivalence of parity $\p-\p'$) \iff\ there are non-zero scalars $\en_n,\en_{n+1/2},\, \en_{n+1},\, \ldots,\, \en_{n+(l-1)/2}$ such that for all $i$ and $j$ in $\oh \bN$ with $j + \frac{3}{2} \le i \le \oh (l-1)$,
\begin{equation} \label{en}
   b^{\p+2j}_{n+i,\, n+j} (\lu)\, \en_{n+i}\ =\ b^{\p'+2j}_{n+i,\, n+j} (\lup)\, \en_{n+j}.
\end{equation}

Moreover, we have the following weaker sufficient conditions for equivalence of the \sq s.  If their \cs\ contains neither of the pairs~(\ref{jump 4}) (which is always the case for $l \le 8$), then they are equivalent \iff\ there exist $e_{n+i}$ such that~(\ref{en}) holds for those $(i, j)$ with $i - j = \frac{3}{2}$, $2$, or $\frac{5}{2}$.  If it does contain one of the pairs~(\ref{jump 4}), then they are equivalent \iff~(\ref{en}) holds in the preceding cases and also in the single case where $i - j = 4$ and $n+j = \frac{1}{4} (-7 \pm \sqrt{33})$.
\end{prop}

\meno {\em Proof.\/}
We first prove the result assuming $\p'=\p$.  Here the \sq s are equivalent \iff\ there is an invertible endomorphism $\en$ of $\bigoplus_{2i = 0}^{l-1} \F^{(\p+2i) \Pi}_{n+i}$ intertwining the \r s $\pi^\lu$ and $\pi^\lup$.  Regard $\en$ as an $l \times l$ matrix with entries $\en_{n+i,\, n+j}$ mapping $\F^{(\p+2j) \Pi}_{n+j}$ to $\F^{(\p+2i) \Pi}_{n+i}$.  Since the restrictions of $\pi^\lu$ and $\pi^\lup$ to $\ds$ are both equal to the block-diagonal \r\ $\bigoplus_{2i = 0}^{l-1} L^{(p+2i) \Pi}_{n+i}|_\ds$, we find that $\en_{n+i,\, n+j}$ must $\ds$-intertwine the tensor density actions.  Applying Lemma~3.2(c) of \cite{Co09a} and remembering that $n$ is non-resonant, we find that $\en$ must be diagonal: $\en_{n+i,\, n+j}$ is a non-zero scalar $\en_{n+i}$ for $i = j$, and zero otherwise.  Corollary~\ref{pi b's} now gives~(\ref{en}) for all $(i, j)$, and the converse is clear.

The second paragraph is analogous to Proposition~7.15 of \cite{CL13}.  The point is that by Theorem~7.5 of \cite{Co09a}, the 1-cochains $\beta_{n+j,\, n+i}$ appearing in Corollary~\ref{pi b's} are 1-cocycles \iff\ $i-j$ is $\frac{3}{2}$, $2$, or $\frac{5}{2}$, or $i - j = 4$ and $n+j = \frac{1}{4} (-7 \pm \sqrt{33})$.  (Theorem~7.5(e) of \cite{Co09a} does not arise for $n$ resonant.)  If the weaker sufficient condition holds we have the following situation: both $\en \circ \pi^\lu \circ \en^{-1}$ and $\pi^\lup$ are \r s of $\K$ on $\bigoplus_{2i = 0}^{l-1} \F^{(\p+2i) \Pi}_{n+i}$.  Their matrices are lower triangular, with the tensor density actions on the diagonal and multiples of the 1-cochains $\beta_{n+j,\, n+i}$ below the diagonal, and these multiples agree in all the entries where $\beta_{n+j,\, n+i}$ is a cocycle.  By the well-known {\em cup equation\/} (see for example Section~4.2 of \cite{Co09a}), this forces the multiples to agree everywhere, proving that $\en$ intertwines the \r s.

Now suppose $\p' \not= \p$.  The \r\ $\pi^\lu$ acts on $\bigoplus_{2i = 0}^{l-1} \F^{(\p+2i) \Pi}_{n+i}$ by the formula of Corollary~\ref{pi b's}.  Recall that the reversed \r\ $(\pi^\lu)^\Pi$ acts on $\bigoplus_{2i = 0}^{l-1} \F^{(\p'+2i) \Pi}_{n+i}$ by $(\pi^\lu)^\Pi(X) = \pi^\lu(X)$.  Observe that at the vector space level, the definition of $\beta_{m,m'}$ is the same whether it is regarded as $\Hom(\F_m, \F_{m'}^{2(m'-m)\Pi})$- or $\Hom(\F_m^\Pi, \F_{m'}^{2(m'-m+1/2)\Pi})$-valued.  It follows that the lower triangular matrix entries of $(\pi^\lu)^\Pi$ are still given by the formula of Corollary~\ref{pi b's}, even though the diagonal entries have been reversed.  Therefore the proof of the $\p' = \p$ case shows that $(\pi^\lu)^\Pi$ is equivalent to $\pi^\lup$ \iff\ there exist $\en_i$ such that~(\ref{en}) holds.  Since $\ep$ intertwines $\pi^\lu$ and $(\pi^\lu)^\Pi$, the first paragraph follows.  The second is proven similarly.  $\Box$

\meno {\em Proofs of Propositions~\ref{L4 prop}, \ref{L5 prop}, \ref{Lle5 prop}, and Theorem~\ref{L6 thm}.\/}
In each of these results, Proposition~\ref{nonres eqvs} shows that the two \r s are equivalent \iff\ there exist non-zero scalars $\en_{n+i}$ such that
\begin{equation} \label{en ratio}
   \frac{\en_{n+i}}{\en_{n+j}}\ =\ 
   \frac{C^{\p'+2j}_{n+i,\, n+j} (\lup)}{C^{\p+2j}_{n+i,\, n+j} (\lu)}\,
   \frac{B^{\p'+2j}_{n+i,\, n+j} (\lup)}{B^{\p+2j}_{n+i,\, n+j} (\lu)}
\end{equation}
for $i-j  = \frac{3}{2}$, $2$, and $\frac{5}{2}$.  Note that only the ratios of the $\en_{n+i}$ enter, and in length~$l$ there are $l-1$ independent ratios: $\en_{n+i}/\en_n$, where $i = \oh, 1, \ldots, \oh(l-1)$.  Clearly the SVC is necessary for the existence of such scalars; assume that it holds.

In Proposition~\ref{L4 prop}, (\ref{en ratio}) need only be solved at $(i,j) = (\frac{3}{2}, 0)$.  This is accomplished simply by setting $\en_{n+3/2} / \en_n$ to the desired value.  In Proposition~\ref{L5 prop} there are three pairs $(i,j)$ at which to solve~(\ref{en ratio}) and four independent ratios, so we can always solve.  Proposition~\ref{Lle5 prop} is a restatement of these two propositions.

In Theorem~\ref{L6 thm} there are six pairs $(i,j)$ and five free ratios.  If $b^{\p+2j}_{n+i,\, n+j}$ vanishes when $(i,j)$ is any of $(\frac{5}{2}, 0)$, $(2, \oh)$, $(\frac{5}{2}, \oh)$, and $(2, 0)$ (in particular, if any of the Pochhammer symbols in the numerators of the $C^{\p+2j}_{n+i,\, n+j}$ vanishes), then~(\ref{en}) is automatically satisfied at that $(i,j)$ and we have enough freedom to solve.  If none of these terms vanishes, we can use five of the six equations~(\ref{en ratio}) to eliminate all scalars $e_\bullet$ from the sixth; the factors $C^{\p+2j}_{n+i,\, n+j}$ then cancel and we obtain the invariant $I_n^\p$.  $\Box$ 

\meno {\em Proofs of Propositions~\ref{SVC prop}, \ref{Lge8 prop}, \ref{Lge15 prop}, and Theorem~\ref{L7 thm}.\/}
For Proposition~\ref{SVC prop}, it is clear from~(\ref{en}) that the SVC is necessary in all lengths.  In Parts~(i) and~(ii) of Theorem~\ref{L7 thm} one of the two modules is ``broken'' by~(\ref{Dsplitting}), causing enough Pochhammer symbols to vanish to reduce the equivalence condition to the SVC and the appropriate conditions from Section~\ref{Lle6}.  For Part~(iii), observe that there are nine equations~(\ref{en ratio}) and six free ratios.  If no $b^\bullet_{\bullet,\bullet}$ vanishes, eliminating the ratios gives the three invariants $I^0_n$, $I^1_{n+1/2}$, and $J^0_n$.  In general, the reader may check that only those $\Simp(\I)$ with no vanishing factors are invariants.

In the first paragraph of Proposition~\ref{Lge8 prop}, if no $b^\bullet_{\bullet,\bullet}$ vanishes then there are $3l - 12$ equations~(\ref{en ratio}).  There are $l-1$ free ratios and therefore $2l - 11$ invariants, which the reader may check are those in~(\ref{arb length}).  If any particular $b^\bullet_{\bullet,\bullet}$ vanishes, then the invariants involving it are replaced by their ratios as in the table following Theorem~\ref{L7 thm}.  No matter how many $b^\bullet_{\bullet,\bullet}$ vanish, if all simplified ratios of the invariants~(\ref{arb length}) that are defined are equal on the two \sq s, then~(\ref{en ratio}) can be solved.  Since these invariants already arise in the length~$7$ ``sub-\sq s'' of the \sq s, the result is proven.

For the second paragraph, check that the equivalence of the two special length~$9$ sub-\sq s forces~(\ref{en}) to hold at $i-j=4$ and $n+j = \frac{1}{4} (-7 \pm \sqrt{33})$.  

In order to prove Proposition~\ref{Lge15 prop}, suppose that the two \sq s are equivalent.  By~(\ref{odd length}), it suffices to treat the case $\p = 0$.  In length~$15$, (\ref{arb length}) gives three groups of five consecutive invariants: $I^0_{n+i}$, $I^1_{n+i+1/2}$, and $J^0_{n+i}$, for $i = 0, 1, 2, 3, 4$.  If $F$ is any one of these invariants, clear denominators in the equation $F(\gd) = F(\gdp)$.  Since the SVC is necessary for equivalence, the resulting equation holds even if some of the relevant $b^\bullet_{\bullet,\bullet}$ are zero.  Applying $n$-difference operators repeatedly in each of the three groups of equations gives $(\gdp) = (\gd)$ (this is an easy computation using software).  $\Box$

\medbreak
Note that if some $b^\bullet_{\bullet,\bullet}$ does vanish, then $(\gd)$ and $(\gdp)$ are restricted to a 1-dimensional curve.  In these cases it should be easy to verify that the only equivalences between \sq s of length $l \ge 7$ are conjugation and the Bol equivalence of Lemma~\ref{Bol}, but we have not found an efficient way to proceed.

\meno {\em Proof of Proposition~\ref{lac 4 prop}.\/}
In the lacunary cases of length $l \le 4$, obvious modifications of the arguments in this section show that the two \sq s are equivalent \iff\ there are non-zero scalars $\en_n, \en_{n+3/2}, \en_{n+2}, \ldots, \en_{n+l/2}, \en_{n+(l+3)/2}$ for which~(\ref{en}) holds where defined.  Clearly the lacunary SVC is necessary for equivalence.

At $l=3$ there are two pairs $(i,j)$ where~(\ref{en}) holds and two free ratios, hence no invariants.  (Indeed, by Proposition~7.8(a) of \cite{Co09a}, in the non-resonant case there is only one indecomposable module composed of $\F^{\p\Pi}_n$, $\F^{(\p+1)\Pi}_{n+3/2}$, and $\F^{\p\Pi}_{n+3}$.)

At $l=4$ there are four pairs and three free ratios, implying one invariant if none of the relevant $b^\bullet_{\bullet,\bullet}$ vanish.  Deduce that the invariant is $M^\p_n$ by eliminating the ratios from the four equations~(\ref{en}).  $\Box$

\section{Remarks on the resonant case}  \label{Resonance}

In summary, we have completely resolved the non-resonant equivalence question in lengths $l \le 6$, but we have not yet found the simplest possible answer in lengths $l \ge 7$.  Moreover, we have not addressed the resonant case: this would require at least in part the computation of the resonant matrix entries of $\pi^\lu$.  As is discussed for $\VR$ in \cite{CL13} and for $\K$ in Section~9 of \cite{Co09a}, these entries may be obtained by taking the appropriate limits of the non-resonant entries.  We conclude with some conjectures concerning the resonant entries and self-dual length~$6$ \sq s.

First we make a conjecture concerning the form of the scalars $b^\p_{m+r, m}(\lu)$ in the non-resonant case: for $\p=0$ or~$1$, $r \in \frac{3}{2} + \oh \bN$, and $m$ non-resonant, we predict that $b^\p_{m+r, m}$ is a product $C^\p_{m+r m} B^\p_{m+r, m}$ as in Theorem~\ref{bCB}, where $C^p_{m+r, m}$ and $B^\p_{m+r, m}$ have the following properties.  For $r \in \frac{3}{2} + \bN$ and $r \in 2 + \bN$, respectively,
\begin{eqnarray*}
   C^\p_{m+r, m} &\propto& \frac{2^{2(r-1)}\, (\delta - m - \p/2)_{\lf r + (1-\p)/2 \rf}}
   {(2m+2r-1)_{r-3/2}\, (2m+r-1/2)\, (2m+r-5/2)_{r-3/2}},
   \\[6pt]
   C^\p_{m+r, m} &\propto& \frac{2^{2(r-1)}\, (\delta - m - \p/2)_{\lf r + (1-\p)/2 \rf}}
   {(2m+2r-1)_{r-1}\, (2m+r-2)_{r-1}},
\end{eqnarray*}
the constants of proportionality depending only on $r$.  The formula for $B^\p_{m+r, m}$ extends~(\ref{Bij}), is polynomial in $\gamma^{1/2}$, $\delta$, and $m$, and is monic in $\gamma^{1/2}$.

It should be easy to verify this prediction from the formula for $b^\p_{m+r, m}$ given in \cite{Co09a}.  This formula is fully simplified only at $r = \frac{3}{2}$, $2$, and~$\frac{5}{2}$, but we expect that the conjecture can be resolved without fully simplifying it for all $r$.

The coefficients $b^\p_{m+r, m}$ have important symmetries arising from Lemmas~\ref{conj} and~\ref{circle dual}, which are given explicitly in Propositions~6.8 and~6.12 of \cite{Co09a}, respectively.  These symmetries appear most clearly when the coefficients are written as functions of $(\gamma^{1/2}, \delta)$:
\begin{eqnarray*}
   b^\p_{m+r, m}(-\gamma^{1/2}, \delta) &=&
   (-1)^{\lf r \rf + 2 \{ r \} \p}\, b^\p_{m+r, m}(\gamma^{1/2}, \delta), \\[6pt]
   b^\p_{m+r, m}(\gamma^{1/2}, -\delta) &=&
   (-1)^{\lf r \rf + 2 \{ r \} \p}\, b^{\p+1+2r}_{-m+1/2, -m-r+1/2}(\gamma^{1/2}, \delta).
\end{eqnarray*}
The $\gamma^{1/2}$-monic polynomials $B^\p_{m+r, m}$ have corresponding symmetries:
\begin{eqnarray*}
   B^\p_{m+r, m}(-\gamma^{1/2}, \delta) &=&
   (-1)^{\lf r \rf + 2 \{ r \} \p}\, B^\p_{m+r, m}(\gamma^{1/2}, \delta), \\[6pt]
   B^\p_{m+r, m}(\gamma^{1/2}, -\delta) &=&
   B^{\p+1+2r}_{-m+1/2, -m-r+1/2}(\gamma^{1/2}, \delta).
\end{eqnarray*}
Note that the $\delta$-symmetry involves a reflection across the {\em antidiagonal\/} of the matrix $\pi^\lu$: the entries $\pi^{\lu, \p}_{m+r, m}$ for which $2m+r = \oh$.

In the resonant case, the restriction of the matrix $\pi^\lu$ to $\ds$ will be ``almost diagonal'': the only non-diagonal entries not zero on $\ds$ will be those on the antidiagonal of the form $\pi^{\lu, \p}_{-m + 1/2,\, m}$ with $m \in -\oh\bN$.  Combining the cup equation with the fact that the non-resonant matrices have zeroes on the first two subdiagonals, it will turn out that in the resonant case the formula for the entries $\pi^{\lu, \p}_{m+r, m}$ on the first and second sub- and super-antidiagonals, where $2m+r \in \{-\oh, 0, 1, \frac{3}{2}\}$, remains the same as in the non-resonant case.  Observe that the denominator of the conjectural formula for $C^\p_{m+r, m}$ vanishes on all resonant entries except these.

By analogy with Theorem~7.10(ii) of \cite{CL13}, we also expect that although the resonant entries $\pi^{\lu, \p}_{-m + 1/2,\, m}$ appearing on the antidiagonal are not multiples of $\beta_{m, -m + 1/2}$, for $m \in -\oh\bZ^+$ they are proportional to $B^\p_{m+r, m}$.  Furthermore, in light of Lemma~\ref{res facs} we predict that the analog of Theorem~7.10(i) is that the $\gamma^{1/2}$-monic polynomial to which the first antidiagonal entry $\pi^{\lu, \p}_{1/2,\, 0}$ is proportional is $\gamma^{\p/2}$, where $\p = 0$ or $1$.  This leads to the following conjecture in the self-dual length~$6$ cases, where $N_6 = 0$, the resonant cases with \cs\
\begin{equation*}
   \bigl\{ \F^{\p\Pi}_{-1},\, \F^{(\p+1)\Pi}_{-1/2},\, \F^{\p\Pi}_0,\,
   \F^{(\p+1)\Pi}_{1/2},\, \F^{\p\Pi}_1,\, \F^{(\p+1)\Pi}_{3/2} \bigr\}.
\end{equation*}

\meno {\bf Conjecture.}
For $\p=0$ or $1$, define $R^\p := \gamma^{\p/2}\, B^\p_{3/2,\, -1}\, /\, B^\p_{3/2,\, 0}\, B^\p_{1/2,\, -1}$.  Then $R^\p$ is an invariant of the \ec\ of the length~$6$ resonant self-dual \sq\ $\SQ^{\delta+1,\, \p,\, 6}_\lu$.  Moreover, $R^\p$ and $I^\p_{-1}$ together with the SVC are complete invariants.

\medbreak
In the case $\p = 0$, one checks easily that
\begin{equation*}
   R^0 = \frac{\gamma - 3/4}{\gamma}, \qquad
   I^0_{-1} = \frac{(\gamma - 3/4)^2}{(\gamma + 1/4)^2 - \delta^2}.
\end{equation*}
Thus when no coefficients vanish, the conjecture predicts that two modules are equivalent \iff\ they have the same values of $\gamma$ and $\delta^2$.

For $\p = 1$ we find
\begin{equation*}
   R^1 = \frac{\gamma (\gamma - 3)}{(\gamma - 3/4)^2 - 9 \delta^2 / 4}, \qquad
   I^1_{-1} = \frac{\gamma(\gamma - 3)}{\gamma^2 - 4\delta^2}.
\end{equation*}
Since $(16/R^1) - (9/I^1_{-1})$ reduces to $(7\gamma - 3)/\gamma$, we obtain the same prediction as in the case $\p = 0$.  Thus in both cases we expect an exceptional equivalence between modules of equal $\gamma$ and opposite $\delta$.

\def\eightit{\it} 
\def\bib{\bf}
\bibliographystyle{amsalpha}

\end{document}